\documentclass[10pt,a4paper]{article}
\usepackage{latexsym,amssymb,amsmath,mathrsfs}
\usepackage{sectsty}
\usepackage{color}
\usepackage{bm}
\usepackage[utf8]{inputenc}
\usepackage[T1]{fontenc}
\usepackage[anchorcolor=blue,%
 bookmarks=true,%
 bookmarksnumbered=true,%
]{hyperref}
\usepackage{cite}
\usepackage{graphicx}

\allsectionsfont{\normalsize\sc\center}

\newtheorem{thm}{Theorem}[section]
\newtheorem{prop}[thm]{Proposition}
\newtheorem{cor}[thm]{Corollary}
\newtheorem{lem}[thm]{Lemma}
\newtheorem{defn}[thm]{Definition}
\newtheorem{rem}[thm]{Remark}
\newtheorem{ex}[thm]{Example}
\newtheorem{assumption}{(A)}

\newenvironment{pf}{\par\begin{trivlist}%
\item[]{\bf Proof.}\ }{\hfill $\square$ \end{trivlist}\par}
\newenvironment{apf}[1]{\par\begin{trivlist}%
\item[]{\bf Proof of #1.}\ }{\hfill $\square$ \end{trivlist}\par}

\makeatletter \@addtoreset{equation}{section} \makeatother

\newcommand{\R}{\mathbb{R}}

\renewcommand{\d}{\mathrm{d}}
\newcommand{\m}{\mathfrak{m} }

\newcommand{\PP}{\bf P}
\newcommand{\EE}{\bf E}

\newcommand{\E}{\mathscr{E}}
\newcommand{\F}{\mathscr{F}}
\newcommand{\loc}{{\rm loc}}
\newcommand{\1}{{\bf 1}}
\newcommand{\eps}{{\varepsilon}}

\newcommand{\<}{\langle}
\renewcommand{\>}{\rangle}

\title{\large\bf $L^p$-Green-tight measures of $L^p$-Kato class for symmetric Markov processes}
\author{Kazuhiro Kuwae\thanks{Department of Applied Mathematics, Fukuoka University,
Fukuoka 814-0180, Japan ({\sf kuwae@}{\sf fukuoka-u.ac.jp}). Supported in part by JSPS Grant-in-Aid for Scientific Research (KAKENHI) 17H02846 and by fund (No.:185001) from the Central Research Institute of Fukuoka University.}
\ \ and\ \ 
Takahiro Mori\thanks{Research Institute for Mathematical Sciences, Kyoto University, Kyoto, 606--8502, 
Japan ({\sf tmori@kurims.kyoto-u.ac.jp}). 
Supported by JSPS Grant-in-Aid for Scientific Research (KAKENHI) 18J21141. 
}
}
\date{}

\begin{document}

\maketitle

\begin{abstract}
In this paper, we introduce the notion of
$L^p$-Green-tight measures of
$L^p$-Kato class
 in the framework of symmetric Markov processes.
The class of $L^p$-Green-tight measures of
$L^p$-Kato class is defined by
the $p$-th power of resolvent kernels.
We first prove that under the $L^p$-Green tightness of the
 measure $\mu$, the embedding of extended Dirichlet space into $L^{2p}(E;\mu)$ is compact under the absolute continuity condition for transient Markov processes, which is an extension of recent seminal work by Takeda.
 Secondly, we prove the coincidence between two classes of
 $L^p$-Green-tightness, one is originally introduced by Zhao, and another one is invented by Chen.
 Finally, we prove that our class of
 $L^p$-Green-tight measures of
$L^p$-Kato class coincides with the class of $L^p$-Green tight measures of Kato class in terms of
Green kernel under the global heat kernel estimates.
We apply our results to $d$-dimensional Brownian motion and
rotationally symmetric relativistic $\alpha$-stable processes on $\R^d$.
\end{abstract}



\section{Introduction}\label{sec:introduction}
The notion of Green-tightness for Kato class potential was
introduced by Zhao~\cite{Zhao:Subcri}
to consider the gaugeability for Feynman-Kac functionals
and the subcriticality of Schr\"odinger operator
$-\frac12\Delta+V$ in the framework of
$d$-dimensional Brownian motion with $d\geq3$.
Before
Zhao~\cite{Zhao:Subcri}, the gaugeability of
Feynman-Kac functionals with Kato class potential
has been considered for absorbing Brownian motions on
bounded open domains (see Zhao~\cite{Zhao:1983}).
Zhao~\cite{Zhao:1989}
also clarified that Kato class potential for
absorbing Brownian motion on a bounded open domain
satisfies the Green-tightness condition in terms of
the Green function of Dirichlet Laplacian on the  domain.
This
was a motivation to formulate the notion of
Green-tight measures of Kato class for
transient symmetric Markov processes.
However,
the Green-tightness as introduced by Zhao~\cite{Zhao:Subcri}
was not enough to develop the theory of gaugeability of
Feynman-Kac functionals and subcriticality of
Schr\"odinger operators for symmetric Markov processes.
To overcome this difficulty,
Chen-Song~\cite{CSgauge2002,CSgauge2003} gave
a new notion of Green-tight smooth measures of Kato class
in the strict sense in the framework of
general $\m$-symmetric transient Markov processes
${\bf X}=(\Omega, X_t, \mathbb{P}_x)$
on a locally compact separable metric space $E$
having a positive Radon measure $\m$ with full support
satisfying the absolute continuity condition with respect to
$\m$.
Moreover,
in Chen~\cite{Chen:gaugeability2002}, this new notion was
refined with remaining value for the gaugeability of
Feynman-Kac functionals and the subcriticality of
Schr\"odinger operators.
Here
${\bf X}$ is said to possess the
\emph{absolute continuity condition} with respect to $\m$
({\bf (AC)} in short) if for any Borel set $B$,
$\m(B)=0$ implies $P_t(x,B)=\mathbb{P}_x(X_t\in B)=0$
for all $t>0$ and $x\in E$.

The refined new class introduced in
\cite{Chen:gaugeability2002} coincides with the class
similarly as defined in \cite{Zhao:Subcri}
not only for $d$-dimensional Brownian motions with $d\geq3$
but also rotationally symmetric $\alpha$-stable processes
with $d>\alpha$.
Chen~\cite{Chen:gaugeability2002}
remarked that if the underlying measure $\m$ of
symmetric Markov process belongs to the Green-tight measures
of Kato class in the original sense of Zhao,
then
it belongs to the new class provided the given process
possesses the strong Feller property
(see \cite[Theorem~4.2]{Chen:gaugeability2002}).
On the other hand,
Kim-Kuwae~\cite[Lemma~4.1]{KK:AnalChara} proved that
the both classes coincide provided the given symmetric
Markov process ${\bf X}=(\Omega, X_t, \mathbb{P}_x)$
possesses the resolvent strong Feller property.
Here
${\bf X}$ is said to possess the
\emph{resolvent strong Feller property} ({\bf (RSF)} in short)
(resp.~\emph{strong Feller property} ({\bf (SF)} in short))
if $R_{\alpha}(\mathscr{B}_b(E))\subset C_b(E)$
for some/any $\alpha>0$
(resp.~$P_t(\mathscr{B}_b(E))\subset C_b(E)$ for any $t>0$),
where
$
 R_{\alpha}f(x)
=
 \mathbb{E}_x
 \left[
  \int_0^{\infty}e^{-\alpha t}f(X_t)\d t
 \right]
=
 \int_{0}^{\infty} e^{-\alpha t} P_t f(x) \d t
$
and
$P_tf(x)=\mathbb{E}_x[f(X_t)]$ for $f\in \mathscr{B}_b(E)$.
Here
$\mathscr{B}_b(E)$ (resp.~$C_b(E)$) denotes the family of
all bounded Borel measurable (resp.~continuous) functions
on $E$.
It is known that the implication
{\bf (SF)}  $\Longrightarrow$
{\bf (RSF)} $\Longrightarrow$
{\bf (AC)}
holds.

Based on the coincidence of two classes of
Green-tight measures of Kato class under {\bf (RSF)},
Takeda~\cite{Takeda:Compact} proved that
the semi-group $(P_t)_{t\geq0}$ of ${\bf X}$
is a compact operator on $L^2(E;\m)$
provided ${\bf X}$ belongs to the Class {\bf (T)}.
Using
the compactness of $P_t$, Takeda proved that
the embedding $\F\hookrightarrow L^2(E;\m)$ is compact.
Here
${\bf X}$ is said to belong to Class {\bf (T)}
if it satisfies that the underlying measure $\m$ of ${\bf X}$
belongs to the $1$-order Green-tight measures of Kato class
in the sense of \cite{Chen:gaugeability2002}
(denoted by $\m\in S_{C\!K_{\infty}}^1({\bf X}^{(1)})$),
${\bf X}$ is irreducible and it possesses {\bf (RSF)}.
Here
${\bf X}^{(1)}$ denotes the $1$-subprocess of ${\bf X}$
defined by ${\bf X}^{(1)}=(\Omega, X_t, \mathbb{P}_x^{(1)})$
with
$
 \mathbb{P}_x^{(1)}(X_t\in A)
=
 e^{-t}\mathbb{P}_x(X_t\in A)
$
for all $t>0$ and $A\in\mathscr{B}(E)$.
${\bf X}$ is said to be \emph{irreducible} ({\bf (I)} in short)
if $B\in \mathscr{B}(E)$ satisfies $P_t\1_Bu=\1_BP_tu$
for all $u\in L^2(E;\m)\cap \mathscr{B}(E)$ and $t>0$,
then $\m(B)=0$ or $\m(B^c)=0$ holds.
Let
$(\F_e, \E)$ be the extended Dirichlet space of ${\bf X}$.
If
${\bf X}$ is transient and $\mu$ is a Dynkin class measure,
Stollmann-Voigt's inequality \cite{SV;potential} implies
the continuity of embedding $\F_e \hookrightarrow L^2(E;\mu)$.
If
${\bf X}$ is transient and $\mu$ is a $0$-order Green-tight
measure of Kato class in the sense of
\cite{Chen:gaugeability2002}, this embedding is compact.

On the other hand,
the notion of $L^p$-Kato class was proposed in \cite{TM:pKato}
by the second named author to obtain
the several probabilistic properties
on the intersection measures.
In \cite{TM:pKato},
Stollmann-Voigt's inequality \cite{SV;potential} is extended
to the measures $\mu$ of $L^p$-Dynkin class,
and
it implies the continuity of the embedding
$\F_e \hookrightarrow L^{2p}(E;\mu)$.
The notion of
$L^p$-Green-tight measures of $L^p$-Kato class in the sense of
Zhao or Chen is a natural extension of usual Green-tight
measures of Kato class in these senses,
and
it is important to investigate the $L^p$-Green-tight measures
of $L^p$-Kato class to establish the compact embedding of
$\F_e \hookrightarrow L^{2p}(E;\mu)$.
For this,
we study the $L^p$-Kato class measures in \cite{KwM:pKato}
under heat kernel estimates for small time.
More precisely,
in \cite{KwM:pKato}, we prepare two classes of
$L^p$-Kato class measures, one is probabilistically defined
and denoted by $S_K^{\,p}({\bf X})$,
another
is analytically defined and denoted by $K_{\nu,\beta}^{\,p}$,
where the parameters $\nu,\beta$ are appeared
in the heat kernel estimates
(see conditions
{\bf (A)}\ref{asmp:SC},
{\bf (A)}\ref{asmp:Bishop} and
{\bf (A)}\ref{asmp:Phi} below).
We prove
the coincidence $S_K^{\,p}({\bf X})=K_{\nu,\beta}^{\,p}$
and provided several criteria for $L^p$-Kato class measures
in \cite{KwM:pKato}.
These results are
fundamental and useful for probabilistic behavior of
intersection measures.
Now we start to state the results of this paper.

\bigskip

Let $S_{C\!K_{\infty}}^{\,p}({\bf X})$
(resp.~$S_{C\!K_{\infty}}^{\,p}({\bf X}^{(1)})$ )
denote the class of $0$-order (resp.~$1$-order)
$L^p$-Green-tight measures of $L^p$-Kato class
in the sense of Chen
(see Definition~\ref{df:GreenTightChen} below).

\bigskip

Our first main theorem in this paper is the following,
which is a natural extension of
\cite[Theorem~4.8]{Takeda:Compact} and
\cite[Corollary~4.3]{TM:pKato}.

\begin{thm}\label{thm:compactEmbedding}
We have the following:
\begin{enumerate}
\item
Suppose that ${\bf X}$ is transient
and it satisfies {\bf (AC)}.
Assume $\mu\in S_{C\!K_{\infty}}^{\,p}({\bf X})$,
or $\mu\in S_K^{\,p}({\bf X})$ and
$\m\in S_{C\!K_{\infty}}^{\,1}({\bf X})$.
Then
$(\E,\F_e)$ is compactly embedded into $ L^{2p}(E;\mu)$.
\item
Suppose that ${\bf X}$ satisfies {\bf (AC)}.
Assume $\mu\in S_{C\!K_{\infty}}^{\,p}({\bf X}^{(1)})$,
or $\mu\in S_K^{\,p}({\bf X})$ and
$\m\in S_{C\!K_{\infty}}^{\,1}({\bf X}^{(1)})$.
Then
$(\E_1,\F)$ is compactly embedded into $ L^{2p}(E;\mu)$.
\end{enumerate}
\end{thm}
The first part of Theorem~\ref{thm:compactEmbedding}(2) is
an extension of recent seminal work \cite{Takeda:Compact}
by Takeda on the case $p=1$ for symmetric Markov processes
in Class {\bf (T)},
and
the second part of Theorem~\ref{thm:compactEmbedding}(2) is
a slight extension of \cite[Corollary~4.3]{TM:pKato}.

\bigskip

Let $S_{C\!K_1}^{\,p}({\bf X})$ denote the class of
$0$-order $L^p$-semi-Green-tight measures of
extended $L^p$-Kato class in the sense of Chen and
$S_{K_1}^{\,p}({\bf X})$ the class of $0$-order
$L^p$-semi-Green-tight measures of extended $L^p$-Kato class
in the sense of Zhao
(see Definition~\ref{df:GreenTightZhao} below).

\bigskip

Our second theorem is the following:

\begin{thm}\label{thm:Coincidence}
Suppose that ${\bf X}$ is transient and it possesses
{\bf (RSF)}.
Then we have
\begin{enumerate}
\item
$
 S_{K_{\infty}}^{\,p}({\bf X})
=
 S_{C\!K_{\infty}}^{\,p}({\bf X})
$.
\quad
{\rm (2)}
$
 S_{K_1}^{\,p}({\bf X})
\cap
 S_{L\!K}^{\,p}({\bf X})
=
 S_{C\!K_1}^{\,p}({\bf X})
\cap
 S_{L\!K}^{\,p}({\bf X})
$.
\end{enumerate}
\end{thm}
This is an extension of \cite[Lemma 4.1]{KimKuwae:GeneralAnal}
which treats the case $p=1$:

\bigskip

To state the third theorem, we need the following conditions
on heat kernel global estimates:
We consider $\nu,\beta\in]0,+\infty[$ and $t_0\in]0,+\infty]$.

\begin{assumption}[Life time condition]\label{asmp:SC}
${\bf X}$ has the following property that
\begin{equation*}
 \lim_{t\to0}
 \sup_{x\in E}
 {\PP}_x(\zeta\leq t)
=:
 \gamma
\in[0,1[.
\end{equation*}
In particular, if ${\bf X}$ is stochastically complete,
that is, ${\bf X}$ is conservative,
then this condition is satisfied with $\gamma=0$.
\end{assumption}
We fix an increasing positive function $V$ on $]0,+\infty[$.
\begin{assumption}[Bishop type inequality]\label{asmp:Bishop}
Suppose that
$r\mapsto V(r)/r^{\nu}$ is increasing or bounded,
and $\sup\limits_{x\in E}\m(B_r(x))\leq V(r)$ for all $r>0$.
\end{assumption}

\begin{assumption}[Upper and lower estimates of heat kernel]\label{asmp:Phi}
Let $\Phi_i$ $(i=1,2)$ be positive decreasing functions
defined on $[0,+\infty[$
which may depend on $t_0$ if $t_0<+\infty$,
let $\nu, \beta >0$ and assume that $\Phi_2$ satisfies
the following condition $H(\Phi_2)$:
\begin{equation*}
 \int_1^{\infty}\frac{(V(t)\lor t^{\nu})\Phi_2(t)}{t}\d t
<
 +\infty
\end{equation*}
and $(\Phi {\EE}_{\nu,\beta})$:
for any $x,y\in E$, $t\in]0,t_0[$
\begin{align*}
 \frac{1}{t^{\nu/\beta}}
 \Phi_1\left(\frac{d(x,y)}{t^{1/\beta}}\right)
\leq
 p_t(x,y)
\leq
 \frac{1}{t^{\nu/\beta}}
 \Phi_2\left(\frac{d(x,y)}{t^{1/\beta}}\right).
\end{align*}
\end{assumption}

\begin{defn}[Kato class \boldmath$K_{\nu,\beta}^{\,p}$]\label{df:KatoGreen}
{\rm Fix $\nu>0$ and $\beta>0$.
For a positive Borel measure $\mu$ on $E$,
$\mu$ is said to be of \emph{$L^p$-Kato ($p$-Kato in short)
class relative to Green kernel}
(write $\mu\in K_{\nu,\beta}^{\,p}$) if
\begin{align*}
 \lim_{r\to0}
 \sup_{x\in E}
 \int_{d(x,y)<r}
  G(x,y)^p
 \mu(\d y)
=0
\ &\text{ for  }\nu\geq\beta,
\\
 \sup_{x\in E}
 \int_{d(x,y)\leq 1}\mu(\d y)
<
 +\infty
\ &\text{ for }\nu<\beta,
\end{align*}
where $G(x,y) := G(d(x,y))$ with
\begin{align*}
 G(r)
:=
 \left\{
  \begin{array}{ll}
   r^{\beta - \nu} & \nu>\beta, \quad r\in]0,+\infty[,
  \\
   \log(r^{-1})    & \nu=\beta, \quad r\in]0,1[.
  \end{array}
 \right.
\end{align*}
For a positive Borel measure $\mu$ on $E$,
$\mu$ is said to be of \emph{local $L^p$-Kato
($p$-Kato in short) class relative to Green kernel}
(write $\mu\in K_{\nu,\beta}^{\,p,\loc}$) if
$\1_G\mu\in K_{\nu,\beta}^{\,p}$
for any relatively compact open set $G$.
Clearly
$
 K_{\nu,\beta}^{\,p}
\subset
 K_{\nu,\beta}^{\,p,\loc}
$.
When $p=1$,
we write $K_{\nu,\beta}$ (resp.~$K_{\nu,\beta}^{\loc}$)
instead of
$K_{\nu,\beta}^{\,1}$ (resp.~$K_{\nu,\beta}^{\,1,\loc}$).
It is shown
in \cite[Lemma~3.6]{KwM:pKato} that
any $\mu\in K_{\nu,\beta}^{\loc}$,
hence any $\mu\in K_{\nu,\beta}$, is a Radon measure.
}
\end{defn}

\begin{defn}[\boldmath$L^p$-Green-tight measures of \boldmath$L^p$-Kato class, \boldmath$K_{\nu,\beta}^{\,p,\infty}$]
{\rm
Assume $\nu>\beta$
and the upper estimate in {\bf (A)}~\ref{asmp:Phi} holds
with $t_0 = +\infty$.
In this case ${\bf X}$ is transient
(see Section~\ref{sec:Equals} below).
A positive Borel measure $\mu$ is said to be a
{\it $L^p$-Green-tight measure of $L^p$-Kato class
(in terms of Green kernel)}
if $\mu\in K_{\nu,\beta}^{\,p}$ and
 for some $o\in E$ it holds that
\begin{align*}
 \lim_{R\to\infty}
 \sup_{x\in E}
 \int_{\{d(y,o)\geq R\}}
  G(x,y)^{\,p}
 \mu(\d y)
=0,
\end{align*}
where $G(x,y)=1/d(x,y)^{\nu-\beta}$.
Denote by $K_{\nu,\beta}^{\,p,\infty}$
the class of all $L^p$-Green-tight measures of
$L^p$-Kato class in terms of Green kernel.
}
\end{defn}

Now we can state the fourth and fifth theorems.

\begin{thm}\label{thm:Equals}
Suppose that
{\bf (A)}\ref{asmp:SC},
{\bf (A)}\ref{asmp:Bishop} and
{\bf (A)}\ref{asmp:Phi} hold.
We also assume $\nu>\beta$ and
the upper estimate in {\bf (A)}\ref{asmp:Phi} holds
with $t_0 = +\infty$.
Then we have the following:
\begin{enumerate}
\item
Assume $\mu(E)<+\infty$,
or $\mu\in S_{K_{\infty}}^1({\bf X})$ with $p>1$.
Then
$\mu\in S_K^{\,p}({\bf X})$ is equivalent to
$\mu\in S_{C\!K_{\infty}}^{\,p}({\bf X})$.
\item
It holds that
$
 S_{C\!K_{\infty}}^{\,p}({\bf X})
=
 S_{K_{\infty}}^{\,p}({\bf X})
$.
If further
the lower estimate in {\bf (A)}\ref{asmp:Phi} also holds
with $t_0 = +\infty$,
then
$
 S_{C\!K_{\infty}}^{\,p}({\bf X})
=
 S_{K_{\infty}}^{\,p}({\bf X})
=
 K_{\nu,\beta}^{\,p,\infty}
$.
\end{enumerate}
\end{thm}

\noindent
Since the $1$-subprocess ${\bf X}^{(1)}$ is always transient,
a similar result also holds.
\begin{thm}\label{thm:Equals(1)}
Suppose that
{\bf (A)}\ref{asmp:SC},
{\bf (A)}\ref{asmp:Bishop} and
{\bf (A)}\ref{asmp:Phi} hold.
Then we have the following:
\begin{enumerate}
\item
Assume $\mu(E)<+\infty$,
or
$\mu\in S_{K_{\infty}}^1({\bf X}^{(1)})$ with $p>1$.
Then
$\mu\in S_K^{\,p}({\bf X})$ is equivalent to
$\mu\in S_{C\!K_{\infty}}^{\,p}({\bf X}^{(1)})$.
\item
It holds that
$
 S_{C\!K_{\infty}}^{\,p}({\bf X}^{(1)})
=
 S_{K_{\infty}}^{\,p}({\bf X}^{(1)})
$.
\end{enumerate}
\end{thm}

The constitution of this paper is as follows.
In Section~\ref{sec:asmp}, we prepare our framework
and explain several definitions.
In Section~\ref{sec:proofCptEmb}, we give the proof of
Theorem~\ref{thm:compactEmbedding}.
In Section~\ref{sec:Coincidence},  we give the proof of
Theorem~\ref{thm:Coincidence}.
In Section~\ref{sec:Equals}, we give the proofs of
Theorems~\ref{thm:Equals} and \ref{thm:Equals(1)}.
In Section~\ref{sec:Examples}, we investigate the various type
of compact embeddings of Sobolev spaces in the framework of
$d$-dimensional Brownian motion and
rotationally symmetric relativistic $\alpha$-stable processes
on $\R^d$.

\section{Preliminary}\label{sec:asmp}
For real numbers $a,b\in\R$, we set $a\lor b := \max\{a,b\}$ and
$a\land b := \min\{a,b\}$.
Let
$(E,d)$ be a locally compact separable metric space and
$\m$ a positive Radon measure with full support.
Let
$E_{\partial} := E\cup\{\partial\}$ be
 the one-point compactification of $E$.
For each $x\in E$ and $r>0$,
denote by $B_r(x) := \{y\in E\mid d(x,y)<r\}$
the open ball with center $x$ and radius $r$.
We
consider and fix a symmetric regular Dirichlet form
$(\E,\F)$ on $L^2(E;\m)$.
Then
there exists a Hunt process
${\bf X}=(\Omega,X_t,\zeta,{\PP}_x)$ such that
for each Borel $u\in L^2(E;\m)$,
$T_tu(x)={\EE}_x[u(X_t)]$ $\m$-a.e.~$x\in E$ for all $t>0$,
where $(T_t)_{t>0}$ is the semigroup associated with $(\E,\F)$.
Here
$\zeta := \inf\{t\geq0\mid X_t=\partial\}$ denotes
the life time of ${\bf X}$.
For a Borel set $B$,
we denote $\sigma_B := \inf\{t>0\mid X_t\in B\}$
(resp.~$\tau_B := \inf\{t>0\mid X_t\notin B\}$)
the \emph{first hitting time to $B$}
(resp.~\emph{first exit time from $B$}).
Throughout this paper, we assume that
${\bf X}$ satisfies {\bf (AC)}.
Under {\bf (AC)},
there exists a jointly measurable function $p_t(x,y)$
defined for all $(t,x,y)\in ]0,+\infty[\times E\times E$
such that
$
 {\EE}_x[u(X_t)]
=
 P_tu(x)
:=
 \int_Ep_t(x,y)u(y)\m(\d y)
$
for any $x\in E$,
bounded Borel function $u$ and $t>0$
(see \cite[Theorem 2]{Yan:1988}).
$p_t(x,y)$
is said to be a \emph{semigroup kernel},
or sometimes called a \emph{heat kernel} of ${\bf X}$
on the analogy of heat kernels of diffusions.
Then
$P_t$ can be extended to contractive semigroups
on $L^p(E;\m)$ for $p\geq1$.
The following are well-known:
\begin{enumerate}
\item[(1)]\label{item:kernel1}
$
 \displaystyle
 p_{t+s}(x,y)
=
 \int_E p_{s}(x,z)p_{t}(z,y)\m(\d z)
$
\quad for all \quad  $x,y \in E$ \quad and \quad $t,s > 0$.
\item[(2)]\label{item:kernel2}
$
 P_{t}(x,\d y)
=
 p_{t}(x,y)\m(\d y)
$
\quad for all \quad $ x \in E$ \quad and \quad $t > 0$.
\item[(3)]\label{item:kernel3}
$
 \displaystyle
 \int_E p_{t}(x,y)\m(\d y) \leq 1
$
\quad for all \quad $ x \in E$ \quad and \quad $t > 0$.
\end{enumerate}
We define
$
 R_{\alpha}(x,y)
:=
 \int_0^{\infty}e^{-\alpha t}p_t(x,y)\d t
$,
and
$
 R(x,y)
:=
 R_0(x,y)
=
 \int_0^{\infty}p_t(x,y)\d t
$
for $\alpha>0$, $x,y\in E$.
$R_{\alpha}(x,y)$ (resp.~$R(x,y)$) is called
the \emph{$\alpha$-order resolvent kernel}
(resp.~\emph{$0$-order resolvent kernel}).

\bigskip

Throughout this paper,
we consider a constant $p\in[1,+\infty[$.
\begin{defn}[{\boldmath$L^p$-Kato class \boldmath$S_K^{\,p}({\bf X})$,
extended \boldmath$L^p$-Kato class \boldmath$S_{E\!K}^{\,p}({\bf X})$,
\boldmath$L^p$-Dynkin class \boldmath$S_D^{\,p}({\bf X})$}]
{\rm A positive Radon measure $\nu$ on $E$
is said to be of \emph{$L^p$-Kato class}
(write $\nu\in S_K^{\,p}({\bf X})$) if
\begin{align}
 \lim_{\alpha\to\infty}
 \sup_{x\in E}
 \int_E
  R_{\alpha}(x,y)^p
 \nu(\d y)=0.\label{SK1}
\end{align}
A positive Radon measure $\nu$ on $E$
is said to be of \emph{extended $L^p$-Kato class}
(write $\nu\in S_{E\!K}^{\,p}({\bf X})$) if
\begin{align}
 \lim_{\alpha\to\infty}
 \sup_{x\in E}
 \int_ER_{\alpha}(x,y)^p\nu(\d y)<1.
\label{SEK1}
\end{align}
A positive Radon measure $\nu$ on $E$
is said to be of \emph{$L^p$-Dynkin class}
(write $\nu\in S_D^{\,p}$) if
\begin{align}
 \sup_{x\in E}\int_ER_{\alpha}(x,y)^p\nu(\d y)
<
 +\infty
\quad \text{ for some }\quad \alpha>0.
\label{SD1}
\end{align}
We denote
$
 R_{\alpha}^{\,p}\nu(x)
:=
 \int_ER_{\alpha}(x,y)^{\,p}\nu(\d y)
$.
Then we see $\nu\in S_K^{\,p}({\bf X})$
(resp.~$\nu\in S_{E\!K}^{\,p}({\bf X})$) if and only if
$
 \lim_{\alpha\to\infty}
 \|R_{\alpha}^{\,p}\nu\|_{\infty}
=
 0
$
(resp.~
$
 \lim_{\alpha\to\infty}
 \|R_{\alpha}^{\,p}\nu\|_{\infty}
<
 1
$),
and $\nu\in S_D^{\,p}({\bf X})$ if and only if
$\|R_{\alpha}^{\,p}\nu\|_{\infty}<+\infty$ for some/all
$\alpha>0$ (see \cite[Proposition~2.6]{TM:pKato}).
Clearly,
$
 S_K^{\,p}({\bf X})
\subset
 S_{E\!K}^{\,p}({\bf X})
\subset
 S_D^{\,p}({\bf X})
$.
A positive Radon measure $\nu$ on $E$
is said to be of \emph{local $L^p$-Kato class}
(write $\nu\in S_{L\!K}^{\,p}({\bf X})$) if
$\1_G\nu\in S_K^{\,p}({\bf X})$
for any relatively compact open set $G$.
When $p=1$,
we may write $S_D({\bf X})$
(resp.
$S_K({\bf X})$,
$S_{E\!K}({\bf X})$,
$S_{L\!K}({\bf X})$) instead of
$S_D^{\,p}({\bf X})$
(resp.
$S_K^{\,p}({\bf X})$,
$S_{E\!K}^{\,p}({\bf X})$,
$S_{L\!K}^{\,p}({\bf X})$) for simplicity.
}
\end{defn}

To the end of Section~\ref{sec:Coincidence},
we basically assume that ${\bf X}$ is transient.

\begin{defn}[\boldmath$L^p$-Green-bounded measures]\label{df:GreenBdd}
{\rm A (positive) Radon measure $\nu$ on $E$ is said to be
of {\it $L^p$-Green-bounded} if
\begin{align}
 \|R^{\,p}\nu\|_{\infty}
:=
 \sup_{x\in E}
 R^{\,p}\nu(x)
<
 +\infty,\label{eq:LpGreenBdd}
\end{align}
where
\begin{align}
 R^{\,p}\nu(x)
:=
 \int_ER(x,y)^{\,p}\nu(\d y).\label{eq:LpGreenPoten}
\end{align}
We denote by $S_{D_0}^{\,p}({\bf X})$ the class of
$L^p$-Green-bounded measures.
}
\end{defn}

We now introduce the notion of $L^p$-Green-tight measures of
$L^p$-Kato class.

\begin{defn}[\boldmath$L^p$-Green-tight measures of \boldmath$L^p$-Kato class in the sense of Zhao]\label{df:GreenTightZhao}
{\rm
A (positive) Radon measure $\nu$ on $E$ is said to be an
{\it $L^p$-Green-tight measure of $L^p$-Kato class
in the sense of Zhao}
if
$
 \nu
\in
 S_{D_0}^{\,p}({\bf X})
\cap
 S_K^{\,p}({\bf X})
$
and for any $\eps>0$ there exists a compact set $K$ such that
\begin{align}
 \|R^{\,p}\1_{K^c}\nu\|_{\infty}
:=
 \sup_{x\in E}R^{\,p}\1_{K^c}\nu(x)
<
 \eps.
\label{eq:LpGreenTightZhao}
\end{align}
A (positive) Radon measure $\nu$ on $E$ is said to be an
{\it $L^p$-semi-Green-tight measure of extended $L^p$-Kato
in the sense of Zhao}
if
$
 \nu
\in
 S_{D_0}^{\,p}({\bf X})
\cap
 S_{E\!K}^{\,p}({\bf X})
$
and there exists a compact set $K$ such that
\begin{align}
 \|R^{\,p}\1_{K^c}\nu\|_{\infty}
:=
 \sup_{x\in E}R^{\,p}\1_{K^c}\nu(x)
<
 1.
\label{eq:semiLpGreenTightZhao}
\end{align}
We denote by $S_{K_{\infty}}^{\,p}({\bf X})$ the class of
$L^p$-Green-tight measures of $L^p$-Kato class
in the sense of Zhao,
and
by $S_{K_1}^{\,p}({\bf X})$ the class of
$L^p$-semi-Green-tight measures of extended $L^p$-Kato class
in the sense of Zhao.
Clearly, we see
$
 S_{K_{\infty}}^{\,p}({\bf X})
\subset
 S_{K_1}^{\,p}({\bf X})
\subset
 S_{D_0}^{\,p}({\bf X})
\cap
 S_{E\!K}^{\,p}({\bf X})
$.
}
\end{defn}

\begin{defn}[\boldmath$L^p$-Green-tight measures of \boldmath$L^p$-Kato class in the sense of Chen]\label{df:GreenTightChen}
{\rm
A (positive) Radon measure $\nu$ on $E$ is said to be an
{\it $L^p$-Green-tight measure of $L^p$-Kato class
in the sense of Chen}
if
for any $\eps>0$ there exist a Borel set $K=K(\eps)$
with $\nu(K)<+\infty$ and $\delta=\delta(\eps)>0$ such that
\begin{align}
 \|R^{\,p}\1_{K^c\cup B}\nu\|_{\infty}
:=
 \sup_{x\in E}R^{\,p}\1_{K^c\cup B}\nu(x)
<
 \eps
\label{eq:LpGreenTightChen}
\end{align}
holds for any Borel subset $B$ of $K$
with $\nu(B)<\delta$.

A (positive) Radon measure $\nu$ on $E$ is said to be an
{\it $L^p$-semi-Green-tight measure of extended $L^p$-Kato
in the sense of Chen}
if
there exist a Borel set $K$ with $\nu(K)<+\infty$
and $\delta>0$ such that
\begin{align}
 \|R^{\,p}\1_{K^c\cup B}\nu\|_{\infty}
:=
 \sup_{x\in E}R^{\,p}\1_{K^c\cup B}\nu(x)
<
 1
\label{eq:semiLpGreenTightChen}
\end{align}
holds for any Borel subset $B$ of $K$
with $\nu(B)<\delta$.
We denote
by $S_{C\!K_{\infty}}^{\,p}({\bf X})$ the class of
$L^p$-Green-tight measures of $L^p$-Kato class
in the sense of Chen,
and by $S_{C\!K_1}^{\,p}({\bf X})$ the class of
$L^p$-semi-Green-tight measures of extended $L^p$-Kato class
in the sense of Chen.
Clearly,
we see
$
 S_{C\!K_{\infty}}^{\,p}({\bf X})
\subset
 S_{C\!K_1}^{\,p}({\bf X})
$.
}
\end{defn}

\begin{rem}
{\rm
In the definitions for $S_{C\!K_{\infty}}^{\,p}({\bf X})$ and
$S_{C\!K_1}^{\,p}({\bf X})$, the Borel set $K$ can be taken
to be a closed (or open), and a compact set.
This is remarked
in \cite[remark after Definition~2.2]{Chen:gaugeability2002}
in the case of $p=1$.
The same thing also holds for general $p$.
}
\end{rem}

Now we state several propositions on
$S_{C\!K_{\infty}}^{\,p}({\bf X})$ and
$S_{C\!K_1}^{\,p}({\bf X})$.

\begin{prop}\label{prop:SemiGreenTightGreenBdd}
Suppose that ${\bf X}$ is transient.
Then it holds that
$
 S_{C\!K_1}^{\,p}({\bf X})
\subset
 S_{D_0}^{\,p}({\bf X})
$.
\end{prop}

\begin{pf}
Suppose $\nu\in S_{C\!K_1}^{\,p}({\bf X})$.
Then
there exist a Borel set $K$ with $\nu(K)<+\infty$ and
$\delta>0$ such that \eqref{eq:semiLpGreenTightChen} holds
for any Borel subset $B$ of $K$ with
$\nu(B)<\delta$.
The Borel set $K$ can be taken to be compact so that
$K$ can be covered by finitely many Borel
(relatively open in $K$) subsets
$\{B_i\}_{i=1}^{\ell}$ of $K$ with $\nu(B_i)<\delta$.
Therefore
\begin{align*}
 \sup_{x\in E}
 \int_E
  R(x,y)^{\,p}
 \nu(\d y)
&\leq
 \sum_{i=1}^{\ell}
 \sup_{x\in E}\int_{K^c \cup B_i}
  R(x,y)^{\,p}
 \nu(\d y)
< \ell 
<
 +\infty.
\end{align*}
\end{pf}

To prove the inclusions
$
 S_{C\!K_{\infty}}^{\,p}({\bf X})
\subset
 S_K^{\,p}({\bf X})
$
and
$
 S_{C\!K_{1}}^{\,p}({\bf X})
\subset
 S_{E\!K}^{\,p}({\bf X})
$,
we need the following lemmas.

\begin{lem}\label{lem:flsc}
Let $f$ be a (nearly) Borel function on $E$.
Suppose that for any $x\in E$ it holds that
\begin{align*}
 {\PP}_x
 \left(
  f(X_0)
 \leq
  \varliminf_{t\to0}f(X_t)
 \right)
=
 1.
\end{align*}
Then $f$ is finely lower semi-continuous on $E$.
\end{lem}
\begin{pf}
Set $B := \{x\in E\mid f(x)>\beta\}$ for $\beta\in\R$.
Then $B$ is a nearly Borel set.
On the event $\{\sigma_{E\setminus B} = 0\}$,
there exists a decreasing sequence $\{t_n\}$ converging to $0$
such that
$X_{t_n}\in E\setminus B$ for all $n\in \mathbb{N}$,
and then
$
 \varliminf_{t\downarrow 0}
 f(X_t)
\leq
 \varliminf_{n\to\infty}
 f(X_{t_n})
\leq
 \beta
$.
By combining this with the assumption,
we have for any $x\in B$
\begin{align*}
 {\PP}_x(\sigma_{E\setminus B} = 0)
\leq
 {\PP}_x(f(X_0)\leq \beta)
=
 {\PP}_x(x\in E\setminus B)
=
 0.
\end{align*}
Thus we have
${\PP}_x(\sigma_{E\setminus B}>0)=1$ for all $x\in B$, i.e.,
$B$ is a finely open set.
Therefore $f$ is finely lower semi-continuous.
\end{pf}

\begin{lem}\label{lem:flscLpPotential}
For any $\nu\in S_D^{\,p}({\bf X})$
(resp.~$\nu\in S_{D_0}^{\,p}({\bf X})$ in transient case),
the $p$-potential function $R_{\alpha}^{\,p}\nu$
(resp.~$R^{\,p}\nu$) is
Borel and finely lower semi-continuous on $E$.
\end{lem}
\begin{pf}
First note that $(x,y)\mapsto R_{\alpha}(x,y)$ is
$\mathscr{B}(E)\times\mathscr{B}(E)$-measurable
in view of the joint measurability of
$(t,x,y)\mapsto p_t(x,y)$.
From this,
we can deduce the Borel measurability of
$x\mapsto \int_ER_{\alpha}(x,y)^{\,p}\nu(\d y)$.
Since
$x\mapsto R_{\alpha}(x,y)$ is a Borel finely continuous
function for a fixed $y\in E$, we have
\begin{align*}
 {\PP}_x
 \left(
  \lim_{t\to0}R_{\alpha}(X_t,y)
 \ne
  R_{\alpha}(X_0,y)
 \right)
=
 0.
\end{align*}
Applying Fubini's theorem to the jointly measurable
function
$
 \1_{\{\lim_{t\to0}R_{\alpha}(X_t,\cdot)\ne R_{\alpha}
 (X_0,\cdot)\}}(\omega,y)
$,
we have
\begin{align*}
 {\PP}_x
 \left(
  \lim_{t\to0}R_{\alpha}(X_t,y)
 =
  R_{\alpha}(X_0,y)
 \quad\nu\text{-a.e.}~y\in E
 \right)
=
 1.
\end{align*}
This implies
\begin{align*}
 \int_E
  R_{\alpha}(X_0,y)^{\,p}
 \nu(\d y)
&=
 \int_E
  \lim_{t\to0}R_{\alpha}(X_t,y)^{\,p}
 \nu(\d y)\\
&\leq
 \varliminf_{t\to0}
 \int_E
  R_{\alpha}(X_t,y)^{\,p}
 \nu(\d y)
\qquad {\PP}_x\text{-a.s.}
\end{align*}
Therefore
$x\mapsto \int_ER_{\alpha}(x,y)^{\,p}\nu(\d y)$
is finely lower semi-continuous
by Lemma~\ref{lem:flsc}.
The proof for the case $\alpha=0$ under the transience of
${\bf X}$ is similar.
\end{pf}

\begin{prop}\label{prop:GreenTightKato}
Suppose that ${\bf X}$ is transient.
Then we have the following inclusions.
\begin{enumerate}
\item[\rm (1)]
$
 S_{C\!K_{\infty}}^{\,p}({\bf X})
\subset
 S_K^{\,p}({\bf X})
$.
\qquad
{\rm (2)}
$
 S_{C\!K_1}^{\,p}({\bf X})
\subset
 S_{E\!K}^{\,p}({\bf X})
$.
\end{enumerate}
\end{prop}
\begin{pf}
\begin{enumerate}
\item
Suppose $\nu\in S_{C\!K_{\infty}}^{\,p}({\bf X})$.
Then for any $\eps>0$, there exist a Borel set $K=K(\eps)$
with $\nu(K)<+\infty$ and $\delta=\delta(\eps)>0$
such that
\eqref{eq:LpGreenTightChen} holds for any subset $B$ of $K$
with $\nu(B)<\delta$.
We may assume that $K$ is a compact set.
Set
$
 B
:=
 \{
  x\in K
 \mid
  \int_ER_{\alpha}(x,y)^{\,p}\nu(\d y)>\eps
 \}
$.
Since
$
 \lim_{\alpha\to\infty}
 \int_ER_{\alpha}(x,y)^{\,p}\nu(\d y)
=
 0
$
for each fixed $x\in E$ and $\nu(K)<+\infty$, we have
$\nu(B)<\delta$ for sufficiently large $\alpha>0$.
Moreover,
the set
$
 K\setminus B
=
 \{
  x\in K
 \mid
  \int_E
   R_{\alpha}(x,y)^{\,p}
  \nu(\d y)
 \leq
  \eps
 \}
$
is a finely closed Borel set
by virtue of Lemma \ref{lem:flscLpPotential}.
Applying Frostman's maximum principle
\begin{align}
 \sup_{x\in E}
 R_{\alpha}^{\,p}\1_{K\setminus B}\nu(x)
=
 \sup_{x\in K\setminus B}
 R_{\alpha}^{\,p}\1_{K\setminus B}\nu(x),
\label{eq:Frostman}
\end{align}
which was proved by \cite[(3.5)]{TM:pKato} for $\alpha=0$,
and its proof remains valid for general $\alpha>0$, we obtain
\begin{align*}
 \sup_{x\in E}
 R_{\alpha}^{\,p}\nu(x)
&\leq
 \sup_{x\in E}
 R_{\alpha}^{\,p}(\1_{K^c\cup B}\nu)(x)
+
 \sup_{x\in E}
 R_{\alpha}^{\,p}(\1_{K\setminus B}\nu)(x)\\
&\leq
 \eps
+
 \sup_{x\in K\setminus B}
 R_{\alpha}^{\,p}(\1_{K\setminus B}\nu)(x)\\
&\leq
 \eps
+
 \sup_{x\in K\setminus B}R_{\alpha}^{\,p}\nu(x)
\leq
 2\eps.
\end{align*}
Hence
$\lim_{\alpha\to\infty}\sup_{x\in E}R_{\alpha}^{\,p}\nu(x)=0$.
\item
Suppose $\nu\in S_{C\!K_1}^{\,p}({\bf X})$.
Then there exist a Borel set $K$ with $\nu(K)<+\infty$
and $\delta>0$ such that
\eqref{eq:semiLpGreenTightChen} holds
for any Borel subset $B$ of $K$
with $\nu(B)<\delta$.
We may assume that $K$ is a compact set.
For
$
 0
<
 \eps
<
 1
-
 \sup_{B\subset K,\nu(B)<\delta}
 \|R_{\alpha}^{\,p}(\1_{K^c\cup B}\nu)\|_{\infty}
$,
we set
\begin{align*}
 B
:=
 \left\{
  x\in K
 \;\left|\;
  \int_E
   R_{\alpha}(x,y)^{\,p}
  \nu(\d y)
 >
  1-\eps
 -
  \sup_{B\subset K,\nu(B)<\delta}
  \|R_{\alpha}^{\,p}(\1_{K^c\cup B}\nu)\|_{\infty}
 \right.
 \right\}.
\end{align*}
Then $\nu(B)<\delta$ for sufficiently large $\alpha>0$ and
$K\setminus B$ is finely closed as proved in (1).
Applying \eqref{eq:Frostman}, we obtain
\begin{align*}
 \sup_{x\in E}
 R_{\alpha}^{\,p}\nu(x)
&\leq
 \sup_{x\in E}
 R_{\alpha}^{\,p}(\1_{K^c\cup B}\nu)(x)
+
 \sup_{x\in E}R_{\alpha}^{\,p}(\1_{K\setminus B}\nu)(x)\\
&\leq
 \sup_{B\subset K,\nu(B)<\delta}
 \|
  R_{\alpha}^{\,p}(\1_{K^c\cup B}\nu)
 \|_{\infty}
+
 \sup_{x\in K\setminus B}
 R_{\alpha}^{\,p}(\1_{K\setminus B}\nu)(x)\\
&\leq
 1-\eps
<
 1.
\end{align*}
Hence
$
 \lim_{\alpha\to\infty}
 \sup_{x\in E}
 R_{\alpha}^{\,p}\nu(x)
<
 1
$.
\end{enumerate}
\end{pf}

\begin{rem}\label{rem:GreenTightInclusion}
{\rm
From Propositions
\ref{prop:SemiGreenTightGreenBdd} and
\ref{prop:GreenTightKato}, we have the following inclusions:
\begin{enumerate}
\item[{\rm (1)}]
$
 S_{C\!K_{\infty}}^{\,p}({\bf X})
\subset
 S_{K_{\infty}}^{\,p}({\bf X})
$.
\qquad (2)
$
 S_{C\!K_1}^{\,p}({\bf X})
\subset
 S_{K_1}^{\,p}({\bf X})
$.
\end{enumerate}
}
\end{rem}

The following proposition is an extension of
\cite[Proposition~4.1]{KK:AnalChara}.

\begin{prop}\label{prop:GreenBddEqui}
Suppose that ${\bf X}$ is transient
and assume $\m\in S_{D_0}^{\,1}({\bf X})$.
Then we have the following:
\begin{enumerate}
\item[{\rm (1)}]
$
 S_{D_0}^{\,p}({\bf X})
=
 S_{D}^{\,p}({\bf X})
$.
\qquad {\rm (2)}\quad
$
 S_{K_{\infty}}^{\,p}({\bf X})
=
 S_{K_{\infty}}^{\,p}({\bf X}^{(1)})
$
\qquad {\rm (3)} \quad
$
 S_{C\!K_{\infty}}^{\,p}({\bf X})
=
 S_{C\!K_{\infty}}^{\,p}({\bf X}^{(1)})
$.
\end{enumerate}
\end{prop}
\begin{pf}
Take a positive Radon measure $\nu$ on $E$
and fix $x\in E$.
By H\"older's inequality, we have
\begin{align*}
 \int_E&
  \left(
   \int_E
    R(x,z)R_{\alpha}(z,y)
   \m(\d z)
  \right)^{p}
 \nu(\d y)\\
&=
 \int_E
  \left\{
   \int_E
    \left(
     \int_E
      R(x,z)R_{\alpha}(z,y)
     \m(\d z)
    \right)^{p-1}
    R_{\alpha}(z,y)
   \nu(\d y)
  \right\}
  R(x,z)
 \m(\d z)\\
&\leq
 \left\{
  \int_E
   \left(
    \int_E
     R(x,z)R_{\alpha}(z,y)
    \m(\d z)
   \right)^{p}
  \nu(\d y)
 \right\}^{\frac{p-1}{p}}
 \|
  R_\alpha ^p \nu
 \|_\infty ^{\frac{1}{p}}
 \|
  R \m
 \|_\infty
,
\end{align*}
that is, it holds that
\begin{align*}
 \left\{
  \int_E
   \left(
    \int_E
     R(x,z)R_{\alpha}(z,y)
    \m(\d z)
   \right)^{p}
  \nu(\d y)
 \right\}^{\frac{1}{p}}
\leq
 \|
  R_\alpha ^p \nu
 \|_\infty ^{\frac{1}{p}}
 \|
  R \m
 \|_\infty
.
\end{align*}
Thanks to the resolvent equation
\begin{align*}
 R(x,y)
=
 R_{\alpha}(x,y)
+
 \alpha
 \int_E
  R(x,z)R_{\alpha}(z,y)
 \m(\d z),
\end{align*}
we have
\begin{align*}
 R^{\,p}\nu(x)^{\frac{1}{p}}
&=
 \left(
  \int_ER(x,y)^{\,p}\nu(\d y)
 \right)^{\frac{1}{p}}\\
&=
 \left(
  \int_E
   \left(
    R_{\alpha}(x,y)
   +
    \alpha
    \int_ER(x,z)R_{\alpha}(z,y)\m(\d z)
   \right)^{p}\nu(\d y)
  \right)^{\frac{1}{p}}\\
&\leq
 \left(
  \int_E R_{\alpha}(x,y)^{p}\nu(\d y)
 \right)^{\frac{1}{p}}
+
 \alpha
 \left(
  \int_E
   \left(
    \int_ER(x,z)R_{\alpha}(z,y)\m(\d z)
   \right)^{p}
  \nu(\d y)
 \right)^{\frac{1}{p}}\\
&\leq
 \|R_{\alpha}^{\,p}\nu\|_{\infty}^{\frac{1}{p}}
+
 \alpha
 \|R_{\alpha}^{\,p}\nu\|_{\infty}^{\frac{1}{p}}
 \|
  R \m
 \|_\infty,
\end{align*}
that is,
$
 \|R^{\,p}\nu\|_{\infty}
\leq
 (1+\alpha\|R \m\|_{\infty})^p
 \|R_{\alpha}^{\,p}\nu\|_{\infty}
$
holds. This implies the each assertion.
\end{pf}

\section{Proof of Theorem~\ref{thm:compactEmbedding}}\label{sec:proofCptEmb}
At the beginning of the section,
we state a useful inequality to estimate the $L^{2p}$-norm
of functions.
The following is essentially proved in
\cite[Theorem~4.1]{TM:pKato}:

\begin{prop}
[\boldmath$p$-version of Stollmann-Voigt's inequality]
\label{prop:p-SV}
Suppose that ${\bf X}$ is transient
and let $\nu\in S_{D_0}^{\,p}({\bf X})$.
Then, it holds that
\begin{equation}
 \|u\|_{L^{2p}(E;\m)}^2
\leq
 \|R^{\,p}\nu\|_{\infty}^{\frac{1}{p}}
 \E(u,u),
\quad\text{ for all }\quad
 u\in\F_e.
\label{eq:StollmannVoigt}
\end{equation}
\end{prop}
In the following,
we omit ``$p$-version'' and simply call this
\emph{Stollmann-Voigt's inequality}.

\begin{lem}\label{lem:Tight}
Set $\mathcal{A}_M := \{u\in\F_e\mid \E(u,u)\leq M\}$.
Suppose $\m\in S_{C\!K_{\infty}}^{\,p}({\bf X})$.
Then it holds that
\begin{align}
 \lim_{L\to\infty}
 \sup_{u\in\mathcal{A}_M}
 \int_{\{u^{2p}\geq L\}}u^{2p}\d\m
=
 0.
\end{align}
\end{lem}
\begin{pf}
Fix $\varepsilon > 0$.
Since
$
 \m
\in
 S_{C\!K_{\infty}}^{\,p}({\bf X})
\subset
 S_{K_{\infty}}^{\,p}({\bf X})
$,
there exist a compact set $K$ and $\delta>0$ such that
\begin{align}
 \sup_{x\in E}
 \int_{K^c\cup B}
  R(x,y)^{\,p}
 \m(\d y)
<
 \varepsilon
\end{align}
holds for any subset $B$ of $K$ with $\m(B)<\delta$.
By applying Stollmann-Voigt's inequality
\eqref{eq:StollmannVoigt},
we see that for sufficiently large $L>0$
\begin{align*}
 \sup_{u\in\mathcal{A}_M}
 \m(\{u^{2p}\geq L\})
&\leq
 \frac{1}{L}
 \sup_{u\in\mathcal{A}_M}
 \int_Eu^{2p}\d\m\\
&\leq
 \frac{1}{L}
 \|R^{\,p} \m\|_{\infty}
 \sup_{u\in\mathcal{A}_M}
 \E(u,u)^p
\leq
 \frac{\|R^{\,p} \m\|_{\infty}M^p}{L}
<
 \delta.
\end{align*}
We regard $B := \{u^{2p}\geq L\}\cap K$
and apply Stollmann-Voigt's inequality
\eqref{eq:StollmannVoigt} to $1_{K^c \cup B} \m$.
Then we have
\begin{align*}
 \int_{\{u^{2p}\geq L\}}u^{2p}\d\m
\leq
 \int_{K^c \cup B}u^{2p}\d\m
\leq
 \|
  R^p \1_{K^c \cup B}
 \|_\infty
 \E(u, u)^p
\leq
 M^p \varepsilon
.
\end{align*}
Hence the conclusion follows.
\end{pf}

\begin{lem}\label{lem:Takeda'}
Suppose $\m\in S_{C\!K_{\infty}}^{\,p}({\bf X})$.
If
$\{g_n\}_{n=1}^{\infty}\subset \F_e$ is a sequence with
$\sup_{n\in\mathbb{N}}\E(g_n,g_n)<+\infty$, and satisfies
$g_n\to g$ $\m$-a.e. (or in $\m$-measure),
then
$g_{n}$ converges to $g$ in $L^{2p}(E;\m)$.
\end{lem}
\begin{pf}
Since $\{g_n\}$ is $\E$-bounded, it is $L^{2p}$-bounded
by Stollmann-Voigt's inequality \eqref{eq:StollmannVoigt}.
The uniform integrability of $\{ g_n ^{2p} \}$
is obtained in Lemma \ref{lem:Tight},
and hence Vitali's Theorem
(see \cite[Theorem~16.6]{Schilling:BM} for example)
gives the conclusion.
\end{pf}

\begin{prop}\label{prop:HKDF}
Let $p_t(x,\cdot)$ be the heat kernel of ${\bf X}$.
Then for each $t>0$ and $x\in E$, $p_t(x,\cdot)\in\F$ and
$\E(p_t(x,\cdot), p_t(x,\cdot))\leq\frac{1}{et}p_t(x,x)$.
Moreover,
if ${\bf X}$ is transient and $\m\in S_{D_0}^{\,p}({\bf X})$
for $p>1$, then $R(p_t(x,\cdot))\in\F_e$ for all $x\in E$.
\end{prop}
\begin{pf}
It is easy to see that $p_t(x,\cdot)\in L^2(E;\m)$ for $t>0$
and $x\in E$, from
\begin{align*}
 \int_E
  p_t(x,y)^2
 \m(\d y)
=
 p_{2t}(x,x)
<
 +\infty.
\end{align*}
Hence $p_t(x,\cdot)=P_{t/2}(p_{t/2}(x,\cdot))\in\F$ and
by \cite[Lemma~4.1]{Takeda:Compact} we see
\begin{align*}
 \E(p_t(x,\cdot),p_t(x,\cdot))
&=
 \E
 (
  P_{t/2}(p_{t/2}(x,\cdot))
 ,
  P_{t/2}(p_{t/2}(x,\cdot))
 )\\
&\leq
 \frac{1}{et}
 (p_{t/2}(x,\cdot),p_{t/2}(x,\cdot))_{\m}
=
 \frac{1}{et}p_t(x,x)<+\infty.
\end{align*}
Suppose further that ${\bf X}$ is transient and
$\m\in S_{D_0}^{\,p}({\bf X})$ for $p>1$.
Since
$\m\in S_{D_0}^{\,p}({\bf X})$ for $p>1$,
by \cite[Theorem~1]{Var85}, \eqref{eq:StollmannVoigt} implies
the ultra-contractivity, i.e.,
there exists $C>0$ such that
\begin{align}
 \|P_t\|_{L^1(E;\m)\to L^{\infty}(E;\m)}
\leq
 Ct^{-\frac{p}{p-1}}
\quad\text{ for \ \  all }\quad t>0,
\label{eq:Ultra1}
\end{align}
equivalently
\begin{align}
 p_t(x,y)
\leq
 Ct^{-\frac{p}{p-1}}
\quad\text{ for \ \ all }\quad x,y\in E
\quad \text{ and }\quad t>0.
\label{eq:Ultra2}
\end{align}
From this, we have
\begin{align}
 \int_E
  R(p_t(x, \cdot))(y)p_t(x,y)
 \m(\d y)
&=
 \int_{2t}^{\infty}p_s(x,x)\d s
\leq
 C(p-1)(2t)^{-\frac{1}{p-1}}
<
 +\infty
\label{eq:FiniteEnergy}
\end{align}
for all $x\in E$. Then
\begin{align*}
 \sup_{\alpha>0}
 \E
 (
  R_{\alpha}(p_t(x,\cdot))
 ,
  R_{\alpha}(p_t(x,\cdot))
 )
&\leq
 \sup_{\alpha>0}
 \E_{\alpha}
 (
  R_{\alpha}(p_t(x,\cdot))
 ,
  R_{\alpha}(p_t(x,\cdot))
 )
=
 \sup_{\alpha>0}
 (p_t(x,\cdot),R_{\alpha}(p_t(x,\cdot)))_{\m}\\
&\leq
 (p_t(x,\cdot),R(p_t(x,\cdot)))_{\m}
<
 +\infty
\end{align*}
so that $R(p_t(x,\cdot))\in\F_e$.
\end{pf}

\begin{lem}\label{lem:compactnessSemigroup}
Suppose $\m\in S_{C\!K_{\infty}}^{\,p}({\bf X})$.
If
$\{g_n\}_{n=1}^{\infty}\subset\F_e(\subset L^{2p}(E;\m))$ is
an $\E$-bounded sequence,
then there exists a subsequence $\{g_{n_k}\}_{k=1}^{\infty}$
such that $\{P_tg_{n_k}\}_{k=1}^{\infty}$
$L^{2p}(E;\m)$-converges.
\end{lem}
\begin{pf}
First note that $P_t(\F_e)\subset \F_e$ and
$\E(P_tg,P_tg)\leq\E(g,g)$ for $g\in\F_e$ holds
by \cite[Lemma~1.5.4]{FOT}.
Suppose
$\m\in S_{C\!K_{\infty}}^{\,p}({\bf X})$ with $p>1$.
In this case $R(p_t(x,\cdot))\in\F_e$ for all $x\in E$ and
$t>0$ by Lemma~\ref{prop:HKDF}.
Since $\{g_n\}\subset \F_e$ is $\E$-bounded,
there exist a subsequence $\{g_{n_k}\}$ and $g\in\F_e$
such that $\{g_{n_k}\}$ $\E$-weakly converges to $g$.
Then
\begin{align*}
 P_tg_{n_k}(x)
&=
 \int_Ep_t(x,y)g_{n_k}(y)\m(\d y)
=
 \E(R(p_t(x,\cdot)),g_{n_k})\\
&\rightarrow
 \E(R(p_t(x,\cdot)),g)
=
 \int_Ep_t(x,y)g(y)\m(\d y)
=
 P_tg(x)
\end{align*}
as $k\to\infty$ for all $x\in E$.
Next
we suppose $p=1$ with
$
 \m
\in
 S_{C\!K_{\infty}}^{\,p}({\bf X})
\subset
 S_{D_0}^{\,p}({\bf X})
=
 S_{D_0}^{\,1}({\bf X})
$.
In this case,
\begin{align*}
 \sup_{n\in\mathbb{N}}
 \|g_n\|_{L^2(E;\m)}^2
\leq
 \|R \m\|_{\infty}
 \sup_{n\in\mathbb{N}}
 \E(g_n,g_n)
<
 +\infty
\end{align*}
implies that there exist a subsequence $\{g_{n_k}\}$ and
$g\in L^2(E;\m)$ such that $\{g_{n_k}\}$ converges to $g$
$L^2(E;\m)$-weakly. Then
\begin{align*}
 P_tg_{n_k}(x)
&=
 \int_Ep_t(x,y)g_{n_k}(y)\m(\d y)
\longrightarrow
 \int_Ep_t(x,y)g(y)\m(\d y)
=
 P_tg(x)
\end{align*}
as $k\to\infty$ for all $x\in E$.
Therefore $P_tg_{n_k}\to P_tg$ in $L^{2p}(E;\m)$
in view of Lemma~\ref{lem:Takeda'}.
\end{pf}

\begin{thm}\label{thm:SemigroupComact}
Suppose that ${\bf X}$ is transient and
$\m\in S_{C\!K_{\infty}}^{\,p}({\bf X})$.
Then
the semigroup $P_t:\F_e\to L^{2p}(E;\m)$ is a compact operator.
\end{thm}
\begin{pf}
Let $\{g_n\}\subset \F_e$ be a sequence $\E$-weakly converges
to $g\in\F_e$.
Then by Lemma~\ref{lem:compactnessSemigroup},
there exists a subsequence $\{g_{n_k}\}_{k=1}^{\infty}$
such that
$\{P_tg_{n_k}\}_{k=1}^{\infty}$ is $L^{2p}(E;\m)$-convergent
to some function $h$.
Since
$\{g_n\}$ $\E$-weakly converges to $g$,
there exists a subsequence $\{g_{n_k}\}$ such that
the Ces\'aro mean $\frac{1}{N}\sum_{k=1}^Ng_{n_k}$
converges to $g$ in $(\F_e,\E)$ by Banach-Saks Theorem,
hence $\frac{1}{N}\sum_{k=1}^NP_tg_{n_k}$ converges to $P_tg$
in $(\F_e,\E)$, in particular, in $L^{2p}(E;\m)$.
Therefore $h=P_tg$ in $L^{2p}(E;\m)$.
\end{pf}

\begin{thm}\label{thm:compactEmbeddingNormal}
We have the following:
\begin{enumerate}
\item
Suppose that ${\bf X}$ is transient and
let $\m\in S_{C\!K_{\infty}}^{\,p}({\bf X})$.
Then
the embedding $\F_e\hookrightarrow L^{2p}(E;\m)$ is compact.
\item
Let $\m\in S_{C\!K_{\infty}}^{\,p}({\bf X}^{(1)})$.
Then the embedding $\F\hookrightarrow L^{2p}(E;\m)$ is compact.
\end{enumerate}
\end{thm}
\begin{pf}
(2) follows from (1). We only prove (1).
Suppose that $\{u_n\}_{n=1}^{\infty}\subset \F_e$
converges to $u\in\F_e$ $\E$-weakly.
We write $u_n^{(k)} := (-k)\lor u_n\land k$ and
$u^{(k)} := (-k)\lor u\land k$ for $k\in\mathbb{N}$.
Then
\begin{align}
 \varlimsup_{n\to\infty}&
 \|u-u_n\|_{L^{2p}(E;\m)}
\nonumber\\
&\leq
 \|u-u^{(k)}\|_{L^{2p}(E;\m)}
+
 \varlimsup_{n\to\infty}
 \|u_n-u_n^{(k)}\|_{L^{2p}(E;\m)}
+
 \varlimsup_{n\to\infty}
 \|u^{(k)}-u_n^{(k)}\|_{L^{2p}(E;\m)}
.
\label{eq:u-u_n}
\end{align}
Regarding the first term of \eqref{eq:u-u_n}, we have
\begin{align}
 \|u-u^{(k)}\|_{L^{2p}(E;\m)} ^{2p}
\leq
 \int_{\{|u|\geq k\}}u^{2p} \d\m
\longrightarrow
 0
\quad\text{ as }\quad k\to\infty.
\label{eq:u-u_n:1}
\end{align}
Regarding the second term of \eqref{eq:u-u_n}, we have
\begin{align}
 \varlimsup_{n\to\infty}
 \|u_n-u_n^{(k)}\|_{L^{2p}(E;\m)} ^{2p}
\leq
 \sup_{n\in\mathbb{N}}
 \int_{\{|u_n|\geq k\}}u_n^{2p} \d\m
\longrightarrow
 0
\quad\text{ as }\quad k\to\infty,
\label{eq:u-u_n:2}
\end{align}
where the convergence follows from the $\E$-boundedness of
$\{u_n\}$ and Lemma \ref{lem:Tight}.
It remains to prove
the convergence of the third term, that is,
$
 \varlimsup_{n\to\infty}
 \|u^{(k)}-u_n^{(k)}\|_{L^{2p}(E;\m)}
\rightarrow
 0
$
as $k \rightarrow \infty$.

Now, we have
\begin{align}
 \|u^{(k)}-P_tu^{(k)}\|_{L^{2p}(E;\m)}
\leq&
 \|u^{(k)}-P_tu^{(k)}\|_{L^{2}(E;\m)}      ^{\frac{1  }{p}}
 \|u^{(k)}-P_tu^{(k)}\|_{L^{\infty}(E;\m)} ^{\frac{p-1}{p}}
\nonumber\\
\leq&
 \left(
  \sqrt{t}
  \E(u^{(k)}, u^{(k)})^{\frac{1}{2}}
 \right)^{\frac{1}{p}}
 \left(
  2k
 \right)^{\frac{p-1}{p}}
\nonumber\\
\leq&
 \left(
  \sqrt{t}
  \E(u, u)^{\frac{1}{2}}
 \right)^{\frac{1}{p}}
 \left(
  2k
 \right)^{\frac{p-1}{p}}
\longrightarrow
 0
\quad\text{ as }\quad t\to0.
\label{eq:uk-Ptuk}
\end{align}
Similarly, we also have
\begin{align}
 \varlimsup_{n\to\infty}
 \|u_n^{(k)}-P_tu_n^{(k)}\|_{L^{2p}(E;\m)}
\leq
 \left(
  \sqrt{t}
  \sup_{n\in\mathbb{N}}
  \E(u, u)^{\frac{1}{2}}
 \right)^{\frac{1}{p}}
 \left(
  2k
 \right)^{\frac{p-1}{p}}
\longrightarrow
 0
\quad\text{ as }\quad t\to0.
\label{eq:unk-Ptunk}
\end{align}
By Theorem \ref{thm:SemigroupComact}, we have
$
 \varlimsup_{n\to\infty}
 \|P_t u-P_t u_n\|_{L^{2p}(E;\m)}
=
 0
$
for each $t>0$, and combining this with
\eqref{eq:u-u_n:1} and
\eqref{eq:u-u_n:2}, we have
\begin{align}
 \varlimsup_{n\to\infty}&
 \|P_tu^{(k)}-P_tu_n^{(k)}\|_{L^{2p}(E;\m)}
\nonumber\\
&\leq
 \varlimsup_{n\to\infty}
 \|P_t u^{(k)}-P_t u\|_{L^{2p}(E;\m)}
+
 \varlimsup_{n\to\infty}
 \|P_t u-P_t u_n\|_{L^{2p}(E;\m)}
+
 \varlimsup_{n\to\infty}
 \|P_tu_n-P_tu_n^{(k)}\|_{L^{2p}(E;\m)}
\nonumber\\
&\leq
 \|u^{(k)}-u\|_{L^{2p}(E;\m)}
+
 \varlimsup_{n\to\infty}
 \|u_n-u_n^{(k)}\|_{L^{2p}(E;\m)}
\longrightarrow
 0
\quad\text{ as }\quad k\to\infty.
\label{eq:Ptuk-Ptunk}
\end{align}
By combining
\eqref{eq:uk-Ptuk},
\eqref{eq:unk-Ptunk} and
\eqref{eq:Ptuk-Ptunk}, we have
\begin{align*}
 \varlimsup_{n\to\infty}&
 \|u^{(k)}-u_n^{(k)}\|_{L^{2p}(E;\m)}
\nonumber
\\
&\leq
 \|u^{(k)}-P_tu^{(k)}\|_{L^{2p}(E;\m)}
+
 \varlimsup_{n\to\infty}
 \|P_tu^{(k)}-P_tu_n^{(k)}\|_{L^{2p}(E;\m)}
+
 \varlimsup_{n\to\infty}
 \|P_tu_n^{(k)}-u^{(k)}\|_{L^{2p}(E;\m)}
,
\end{align*}
which converges to $0$ by letting $t\rightarrow 0$
and then $k\rightarrow \infty$.
Therefore we complete the proof.
\end{pf}

\begin{apf}{Theorem~\ref{thm:compactEmbedding}}
(2) is an easy consequence of (1).
We only prove (1).
Suppose first $\mu\in S_{C\!K_{\infty}}^{\,p}({\bf X})$.
Then $\mu\in S_1({\bf X})$,
that is, $\mu$ is a smooth measure in the strict sense
with respect to ${\bf X}$.
Let
$A_t^{\mu}$ be a positive continuous additive functional 
in the strict sense associated with $\mu$
and $F$ its fine support.
Let
$(\check{\bf X},\mu)$ be the time changed process
 of ${\bf X}$ with respect to $A^\mu$ and
$\check{R}(x,y)$ be its $0$-order resolvent kernel.
It is proved in \cite[proof of Lemma~3.4]{TM:pKato} that
for all $x\in F$
\begin{align}
 \check{R}(x,y)
=
 R(x,y)
\quad\mu\text{-a.e.}~y\in F.
\label{eq:resolventkernelTimeChange}
\end{align}
Then we see
$\mu\in S_{C\!K_{\infty}}^{\,p}(\check{\bf X})$.
Let
$(\check{\E},\check{\F}_e)$ be the extended Dirichlet space
associated to the time changed process $(\check{\bf X},\mu)$
(see \cite[(6.2.7)]{FOT}).
$(\check{\E},\check{\F}_e)$ is given by
\begin{align*}
 \left\{
 \begin{array}{rl}
 \check{\F}_e\hspace{-2mm}
&=
 \{\varphi=u|_Y\;\mu\text{-a.e.}\mid u\in\F_e\}\\
 \check{\E}(\varphi,\varphi)\hspace{-2mm}
&=
 \E(\mathscr{P}u,\mathscr{P}u)\quad
 \varphi\in\check{\F}_e,\;
 \varphi=u|_Y\;\mu\text{-a.e.},\;
 u\in\F_e.
 \end{array}
 \right.
\end{align*}
Here $\mathscr{P}$ denotes the orthogonal projection on
$\mathscr{H}_F := (\F_{e,E\setminus F})^{\perp}$
in the Hilbert space $(\E,\F_e)$,
$
 \mathscr{P}u(x)
:=
 H_F\tilde{u}(x)
=
 {\EE}[\tilde{u}(X_{\sigma_F})]
$
and $Y={\rm supp}[\mu]$ is the topological support of $\mu$.
Here
$
 \F_{e,E\setminus F}
:=
 \{u\in\F_e\mid \tilde{u}=0\text{ q.e.~on }F\}
$.
Now suppose that $\{u_n\}\subset \F_e$
$\E$-weakly converges to $u\in\F_e$ in $(\E,\F_e)$.
Let $\varphi_n,\varphi\in\check{\F}_e$
with $\varphi_n := u_n|_Y$ and $\varphi=u|_Y$ $\mu$-a.e.,
and take any $\psi\in\check{\F}_e$ with $\psi=v|_Y$ $\mu$-a.e.
Thus we have
\begin{align*}
 \check{\E}(\varphi_n-\varphi,\psi)
&=
 \E(\mathscr{P}u_n-\mathscr{P}u,\mathscr{P}v)
=
 \E(u_n-n,\mathscr{P}v)
\longrightarrow
 0
\quad \text{ as }\quad n\to\infty,
\end{align*}
that is, $\{\varphi_n\}$
$\check{\E}$-weakly converges to $\varphi$.
Since
$\check{\F}_e$ is compactly embedded into $L^{2p}(Y;\mu)$
by Theorem~\ref{thm:compactEmbeddingNormal},
we conclude
\begin{align*}
 \int_E
  |u_n-u|^{2p}
 \d\mu
=
 \int_Y
  |\varphi_n-\varphi|^{2p}
 \d\mu
\longrightarrow
 0
\quad\text{ as }\quad n\to\infty.
\end{align*}

Next we suppose $\mu\in S_K^{\,p}({\bf X})$ and
$\m\in S_{C\!K_{\infty}}^{\,1}({\bf X})$.
Then $\m\in S_{D_0}^{\,1}({\bf X})$, hence $\F_e=\F$ by
\begin{align*}
 \|u\|_{L^2(E;\m)}^2
\leq
 \|R1\|_{\infty}
 \E(u,u)
\quad \text{ for }\quad u\in\F_e.
\end{align*}
Applying the former result to this case, we can see
the compact embedding $\F_e\hookrightarrow L^2(E;\m)$.
Moreover,
let $\{u_n\}\subset \F_e=\F$ be an $\E$-bounded sequence.
Then
there exist a subsequence $\{u_{n_k}\}$ and $u\in\F_e=\F$
such that $\{u_{n_k}\}$ $\E$-weakly and
$L^2(E;\m)$-strongly converges to $u$.
Since $\mu\in S_K^{\,p}({\bf X})$, we have
\begin{align*}
 \lim_{k,l\to\infty}
 \|u_{n_k}-u_{n_l}\|_{L^{2p}(E;\mu)}^2
&\leq
 \lim_{k,l\to\infty}
 \|R_{\alpha}^{\,p}\mu\|_{\infty}
 \E_{\alpha}(u_{n_k}-u_{n_l},u_{n_k}-u_{n_l})\\
&\leq
 \|R_{\alpha}^{\,p}\mu\|_{\infty}
 \sup_{n\in\mathbb{N}}
 \E(u_n,u_n)
\longrightarrow
 0
\quad\text{ as }\quad \alpha\to\infty.
\end{align*}
This implies the $L^{2p}(E;\mu)$-strong convergence of
$\{u_{n_k}\}$.
\end{apf}

The next proposition is an addendum.

\begin{prop}\label{prop:CompactEmbedPtCompact}
Suppose that the embedding
$\F\hookrightarrow L^{2p}(E;\m)$ is continuous.
Then the following statements are equivalent.
\begin{enumerate}
\item
The embedding $\F\hookrightarrow L^{2p}(E;\m)$ is compact.
\item
$P_t:\F\to L^{2p}(E;\m)$ is a compact operator for $t>0$.
\item
$P_t:L^2(E;\m)\to L^{2p}(E;\m)$ is a compact operator
for $t>0$.
\end{enumerate}
\end{prop}
\begin{pf}
(3)$\Longrightarrow$(2) is trivial.
The proof of (2)$\Longrightarrow$(1) is already done
in the proof of Theorem~\ref{thm:compactEmbeddingNormal}(1)
based on Theorem~\ref{thm:SemigroupComact}.
We prove (1)$\Longrightarrow$(3) only.
Since
$\E(P_tf,P_tf)\leq\frac{1}{2et}(f,f)_{\m}$,
$t>0$, $f\in L^2(E;\m)$ by \cite[Lemma~4.1]{Takeda:Compact},
$P_t:L^2(E;\m)\to\F$ is a bounded operator.
Hence
$P_t: L^2(E;\m)\to L^{2p}(E;\m)$ is the composition of
these operators so that it is compact.
\end{pf}

\section{Proof of Theorem~\ref{thm:Coincidence}}\label{sec:Coincidence}

\begin{apf}{Theorem~\ref{thm:Coincidence}}
By Remark~\ref{rem:GreenTightInclusion}, we already know
$
 S_{C\!K_{\infty}}^{\,p}({\bf X})
\subset
 S_{K_{\infty}}^{\,p}({\bf X})
$
and
$
 S_{C\!K_1}^{\,p}({\bf X})
\subset
 S_{K_1}^{\,p}({\bf X})
$
by Proposition~\ref{prop:GreenTightKato}.
Suppose that ${\bf X}$ possesses {\bf (RSF)}.

We first prove (1).
Take $\nu\in S_{K_{\infty}}^{\,p}({\bf X})$.
Then
$
 \nu
\in
 S_K^{\,p}({\bf X})
\cap
 S_{D_0}^{\,p}({\bf X})
$
and assume $\nu\notin S_{C\!K_{\infty}}^{\,p}({\bf X})$.
Then there is an $\eps>0$ such that
for any $\delta>0$ and any compact set $K$ with
$\sup_{x\in E}R^{\,p}(\1_{K^c}\nu)(x)<\eps/2$,
there exists a Borel subset $B$ of $K$ with $\nu(B)<\delta$
satisfying
$\sup_{x\in E}R^{\,p}(\1_B\nu)(x)\geq\frac{\eps}{2}$.
Let $\{B_n\}_{n=1}^{\infty}$ be a sequence of
such Borel subsets of $K$ with $\nu(B_n)<1/2^n$.
Define $A_n := \bigcup_{k=n}^{\infty}B_k$.
Then $\nu(A_n)<1/2^{n-1}$.
We have for any $n\in \mathbb{N}$,
\begin{align*}
 \frac{\eps}{2}
&\leq
 \sup_{x\in E}
 \int_{B_n}R(x,y)^{\,p}\nu(\d y)\\
&\leq
 \sup_{x\in E}
 \int_{A_n}R(x,y)^{\,p}\nu(\d y)\\
&=
 \sup_{M>0}
 \sup_{x\in E}
 \int_{A_n}R(x,y)\left(R(x,y)\land M\right)^{\,p-1}\nu(\d y).
\end{align*}
Now set
\begin{align*}
 a_n(M)
:&=
 \sup_{x\in E}
 \int_{A_n}
  R(x,y)
  \left(R(x,y)\land M\right)^{\,p-1}
 \nu(\d y)
 \left(
  \leq M^{p-1}\|R\1_K\nu\|_{\infty}
 \right),\\
 a_{\infty}(M)
:&=
 \lim_{n\to\infty}
 a_n(M),
\end{align*}
and apply Terkelsen's minimax principle
(see \cite[Corollary 1]{Terkelsen:1972})
for the continuous function $n\mapsto a_n(M)$
on the compact set $\mathbb{N}\cup\{\infty\}$.
Then we have
\begin{align*}
 \frac{\eps}{2}
\leq
 \lim_{n\to\infty}
 \sup_{M>0}
 a_n(M)
&=
 \min_{n\in\mathbb{N}\cup\{\infty\}}
 \sup_{M>0}
 a_n(M)\\
\\
&=
 \sup_{M>0}
  \min_{n\in\mathbb{N}\cup\{\infty\}}
  a_n(M)
=
 \sup_{M>0}
 \lim_{n\to\infty}
 a_n(M).
\end{align*}
We will see that
$\lim_{n\to\infty}a_n(M)=0$ for any $M>0$,
which gives a contradiction.
Recall $\nu\in S_K^{\,p}({\bf X})$.
Since $\nu(K)<+\infty$,
we have $\1_K\nu\in S_K^{\,1}({\bf X})$.
The function
$
 f_n(x)
:=
 M^{p-1}
 \int_{A_n}
  R(x,y)
 \nu(\d y)
$
is bounded and continuous in view of {\bf (RSF)}
for the time changed process $(\check{\bf X},\1_K\nu)$
(see \cite[Lemma~4.1]{KK:AnalChara}).
Since $\nu(\bigcap_{n=1}^{\infty}A_n)=0$,
we have $f_n(x)\to0$ as $n\to\infty$ for each $x\in E$.
By use of Dini's theorem, we have
$
 \lim_{n\to\infty}
 a_n(M)
\leq
 \lim_{n\to\infty}
 \|f_n\|_{\infty}
=
 0
$,
because
$
 \|f_n\|_{\infty}
=
 \sup_{x\in E}
 f_n(x)
=
 \sup_{x\in K}
 f_n(x)
$
in view of Frostman's maximum principle.
Hence we complete the proof of (1).

\vspace{0.5eM}
We next prove (2).
Take
$
 \nu
\in
 S_{K_1}^{\,p}({\bf X})
\cap
 S_{L\!K}^{\,p}({\bf X})
$,
then
$
 \nu
\in
 S_{D_0}^{\,p}({\bf X})
\cap
 S_{E\!K}^{\,p}({\bf X})
$
by the definition of $S_{K_1}^{\,p}({\bf X})$.
Assume $\nu\notin S_{C\!K_1}^{\,p}({\bf X})$.
Then
for any $\delta>0$ and any compact set $K$ with
$\sup_{x\in E}R^{\,p}(\1_{K^c}\nu)(x)<1$,
there exists a Borel subset $B$ of $K$ with $\nu(B)<\delta$
satisfying
$
 \sup_{x\in E}R^{\,p}(\1_B\nu)(x)
\geq
 1-\sup_{x\in E}R^{\,p}(\1_{K^c}\nu)(x)
$.
Let $\{B_n\}_{n=1}^{\infty}$ be a sequence of
 such Borel subsets of $K$ with $\nu(B_n)<1/2^n$.
Define $A_n := \bigcup_{k=n}^{\infty}B_k$.
Then $\nu(A_n)<1/2^{n-1}$.
In the same way as in the proof of (1),
we have for each $n\in\mathbb{N}$,
\begin{align*}
 1-\sup_{x\in E}R^{\,p}(\1_{K^c}\nu)(x)
&\leq
 \sup_{x\in E}
 \int_{B_n}
  R(x,y)^{\,p}
 \nu(\d y)\\
&=
 \sup_{M>0}
 \sup_{x\in E}
 \int_{A_n}
  R(x,y)
  \left(
   R(x,y)\land M
  \right)^{p-1}
 \nu(\d y).
\end{align*}
Then one can obtain a similar contradiction
as in the proof of (1) by replacing $\eps/2$
with $1-\sup_{x\in K}R^{\,p}(\1_{K^c}\nu)(x)$.
Note that
the time changed process $(\check{\bf X},\1_K\nu)$
possesses {\bf (RSF)}, because $\1_K\nu\in S_K^{\,p}({\bf X})$
implies $\1_K\nu\in S_K({\bf X})$ by $\nu(K)<+\infty$.
\end{apf}

\section{Proofs of  Theorems~\ref{thm:Equals} and \ref{thm:Equals(1)}}\label{sec:Equals}

In this section we prove Theorems \ref{thm:Equals}
and \ref{thm:Equals(1)}.
First note that
under the conditions in Theorem~\ref{thm:Equals},
we have the coincidence
$S_K^{\,p}({\bf X})=K_{\nu,\beta}^{\,p}$
proved in \cite{KwM:pKato}.
Moreover,
by use of \cite[Lemma~4.3]{KwT:Katounderheat},
${\bf X}$ is transient
provided
$\nu>\beta$ and
the upper estimate in {\bf (A)}\ref{asmp:Phi} holds
 with $t_0 = +\infty$.
Indeed,
from \cite[Lemma~4.3]{KwT:Katounderheat},
we easily see that there exists $C>0$ such that
$R(x,y)\leq CG(x,y)\leq +\infty$ for $x,y\in E$.
Since $\m\in S_K^1({\bf X})=K_{\nu,\beta}^1$,
we have
\begin{align*}
 R f (x)
&\leq
 C
 \int_E
  G(x,y)f(y)
 \m(\d y)
\\
&=
 C
 \int_{B_r(x)}
  \frac{f(y)}{d(x,y)^{\nu-\beta}}
 \m (\d y)
+
 C
 \int_{B_r(x)^c}
  \frac{f(y)}{d(x,y)^{\nu-\beta}}
 \m (dy)
\\
&\leq
 C \|f\|_{L^\infty}
 \sup_{x\in E}
 \int_{B_r(x)}
  \frac{\m(\d y)}{d(x,y)^{\nu-\beta}}
+
 C \|f\|_{L^1}
 \frac{1}{r^{\nu-\beta}}
<
 +\infty
\end{align*}
for any $f\in L^1(E; \m) \cap L^\infty(E; \m)$
which are positive $\m$-a.e. on $E$.
This implies the transience of ${\bf X}$
in the sense of \cite[Lemma 1.5.1]{FOT}.

\begin{apf}{Theorem~\ref{thm:Equals} (1)}
Note that
\cite[Lemma~4.3(3)]{KwT:Katofunc}
under these conditions gives that,
there exists $C>0$ such that
$R(x, y) \leq CG(x, y)$
for all $x, y\in E$.
We already know
$
 S_{C\!K_{\infty}}^{\,p}({\bf X})
\subset
 S_K^{\,p}({\bf X})
$
by Proposition \ref{prop:GreenTightKato}.
Assume that
$\mu\in S_K^{\,p}({\bf X})$ is a finite measure.
Then for any $A\in\mathscr{B}(E)$,
\begin{align}
 \int_AG(x,y)^{\,p}\mu(\d y)
&\leq
 \int_{B_r(x)}G(x,y)^{\,p}\mu(\d y)
+
 \int_{A\setminus B_r(x)}G(x,y)^{\,p}\mu(\d y)
\nonumber\\
&\leq
 \sup_{x\in E}
 \int_{B_r(x)}G(x,y)^{\,p}\mu(\d y)
+
 r^{p(\beta-\nu)}\mu(A).
\label{eq:GreenEstimate}
\end{align}
Fix $\eps >0$.
We choose small $r>0$ so that
$\sup_{x\in E}\int_{B_r(x)}G(x,y)^{\,p}\mu(\d y)<\eps/3$,
a compact set $K$ satisfying
$r^{p(\beta-\nu)}\mu(K^c)<\eps/3$ and $\delta>0$ so that
$r^{p(\beta-\nu)}\delta<\eps/3$.
Applying \eqref{eq:GreenEstimate} to $A=K^c\cup B$,
we have for $B\subset K$ with $\mu(B)<\delta$
\begin{align*}
 \sup_{x\in E}
 \int_{K^c\cup B}G(x,y)^{\,p}\mu(\d y)
&\leq
 \sup_{x\in E}
 \int_{B_r(x)}G(x,y)^{\,p}\mu(\d y)
+
 r^{p(\beta-\nu)}
 \mu(K^c\cup B)\\
&\leq
 \frac{\eps}{3}
+
 \frac{\eps}{3}
+
 \frac{\eps}{3}
=
 \eps.
\end{align*}
Thus we have $\mu\in  S_{C\!K_{\infty}}^{\,p}({\bf X})$.

Next suppose $p>1$ and
$\mu\in S_K^{\,p}({\bf X})\cap S_{K_{\infty}}^1({\bf X})$.
Then for any $A\in\mathscr{B}(E)$,
\begin{align}
 \int_AG(x,y)^{\,p}\mu(\d y)
&\leq
 \int_{B_r(x)}G(x,y)^{\,p}\mu(\d y)
+
 \int_{A\setminus B_r(x)}G(x,y)^{\,p}\mu(\d y)
\nonumber\\
&\leq
 \sup_{x\in E}
 \int_{B_r(x)}G(x,y)^{\,p}\mu(\d y)
+
 r^{(p-1)(\beta-\nu)}
 \sup_{x\in E}
 \int_{A\setminus B_r(x)}G(x,y)\mu(\d y).
\label{eq:GreenEstimateGreenTight}
\end{align}
Fix $\eps > 0$.
We choose small $r>0$ so that
$\sup_{x\in E}\int_{B_r(x)}G(x,y)^{\,p}\mu(\d y)<\eps/3$,
a compact set $K$ satisfying
$
 r^{(p-1)(\beta-\nu)}
 \sup_{x\in E}
 \int_{K^c}G(x,y)\mu(\d y)
<
 \eps/3
$,
and $\delta>0$ so that $r^{p(\beta-\nu)}\delta<\eps/3$.
Applying \eqref{eq:GreenEstimateGreenTight} to $A=K^c\cup B$,
we have for $B\subset K$ with $\mu(B)<\delta$
\begin{align*}
 \int_{K^c\cup B}G(x,y)^{\,p}\mu(\d y)
&\leq
 \sup_{x\in E}\int_{B_r(x)}G(x,y)^{\,p}\mu(\d y)
+
 r^{(p-1)(\beta-\nu)}
 \sup_{x\in E}
 \int_{(K^c\cup B)\setminus B_r(x)}G(x,y)\mu(\d y)\\
&\leq
 \sup_{x\in E}
 \int_{B_r(x)}G(x,y)^{\,p}\mu(\d y)
+
 r^{(p-1)(\beta-\nu)}
 \sup_{x\in E}\int_{K^c}G(x,y)\mu(\d y)\\
&\hspace{2cm}+
 r^{p(\beta-\nu)}\mu(B)\\
&\leq
 \frac{\eps}{3}
+
 \frac{\eps}{3}
+
 \frac{\eps}{3}
=
 \eps.
\end{align*}
Thus we have $\mu\in S_{C\!K_{\infty}}^{\,p}({\bf X})$.
\end{apf}

\begin{apf}{Theorem~\ref{thm:Equals} (2)}
We first show the first half of the claim.
By applying \cite[Lemma~4.3(3)]{KwT:Katounderheat} with
$t_0=+\infty$,
there exists $C>0$ such that
$R(x,y)\leq C G(x,y)$ for all $x,y\in E$.
The inclusion
$
 S_{C\!K_{\infty}}^{\,p}({\bf X})
\subset
 S_{K_{\infty}}^{\,p}({\bf X})
$
is given in Remark \ref{rem:GreenTightInclusion}.
It suffices to show the converse inclusion.
Take
$
 \mu
\in
 S_{K_{\infty}}^{\,p}({\bf X})
\subset
 S_K^{\,p}({\bf X})
$.
Then,
for any $\eps>0$ there exists a compact set $K$ such that
$
 \sup_{x\in E}
 \int_{K^c}R(x,y)^p\mu(\d y)
<
 \eps/2
$.
Since  $\mu\in S_K^{\,p}({\bf X})=K_{\nu,\beta}^{\,p}$,
for any $\eps>0$, there exists $r>0$ such that
\begin{align*}
 \sup_{x\in E}
 \int_{B_r(z)}R(x,y)^p\mu(\d y)
\leq
 C^p
 \sup_{x\in E}
 \int_{B_r(x)}G(x,y)^p\mu(\d y)
<
 \frac{\eps}{2}.
\end{align*}
Take a small $\delta>0$ so that $\delta/r^{p(\nu-\beta)}<\eps/2$.
Then, for any Borel subset $B$ of $K$ with $\mu(B)<\delta$
we have
\begin{align*}
 \sup_{x\in E}
 \int_B
  R(x,y)^{\,p}
 \mu(\d y)
&\leq
 C^p
 \sup_{x\in E}
 \int_B
  G(x,y)^{\,p}
 \mu(\d y)
\\
&\leq
 C^p
 \left(
  \sup_{x\in E}
  \int_{B\cap B_r(x)}
   G(x,y)^{\,p}
  \mu(\d y)
 +
  \sup_{x\in E}
  \int_{B\cap B_r(x)^c}
   G(x,y)^{\,p}
  \mu(\d y)
 \right)
\\
&\leq
 C^p
 \left(
  \frac{\eps}{2}
 +
  \frac{\mu(B)}{r^{p(\nu-\beta)}}
 \right)
<
 C^p \eps.
\end{align*}
So we obtain $\mu\in S_{C\!K_{\infty}}^{\,p}({\bf X})$
for this case.

Next we suppose
the lower estimate in {\bf (A)}\ref{asmp:Phi} also holds
with $t_0 = +\infty$.
Then there exist $C_1,C_2>0$ such that
$C_1G(x,y)\leq R(x,y)\leq C_2G(x,y)$ for all $x,y\in E$
by applying \cite[Lemmas~4.1(3) and 4.3(3)]{KwT:Katounderheat} with
$t_0=+\infty$.
The inclusion
$
 S_{C\!K_{\infty}}^{\,p}({\bf X})
\subset
 S_{K_{\infty}}^{\,p}({\bf X})
$
is given in Remark \ref{rem:GreenTightInclusion}.
To prove
$
 S_{K_{\infty}}^{\,p}({\bf X})
\subset
 K_{\nu,\beta}^{\,p,\infty}
$
let $o\in E$ and a compact set $K$.
For sufficiently large $R>0$ with $K\subset B_R(o)$
we have
\begin{align*}
 \sup_{x\in E}
 \int_{B_R(o)^c}G(x,y)^{\,p}\mu(\d y)
\leq
 C_1^{-1}
 \sup_{x\in E}
 \int_{K^c}R(x,y)^{\,p}\mu(\d y)
,
\end{align*}
hence the inclusion
$
 S_{K_{\infty}}^{\,p}({\bf X})
\subset
 K_{\nu,\beta}^{\,p,\infty}
$
holds.
It remains to prove the inclusion
$
 K_{\nu,\beta}^{\,p,\infty}
\subset
 S_{C\!K_{\infty}}^{\,p}({\bf X})
$.
Take $\mu\in K_{\nu,\beta}^{\,p,\infty}$.
Then for any $\eps>0$, there is $R>0$ such that
\begin{align*}
 \sup_{x\in E}
 \int_{B_R(o)^c}G(x,y)^{\,p}\mu(\d y)
<
 \frac{\eps}{2}.
\end{align*}
We choose small $r>0$ so that
$\sup_{x\in E}\int_{B_r(x)}G(x,y)^{\,p}\mu(\d y)<\eps/4$
and $\delta>0$ so that $r^{p(\beta-\nu)}\delta<\eps/4$.
Then applying \eqref{eq:GreenEstimate} to $A=B$, we have
\begin{align*}
 \sup_{B\subset B_R(o),\mu(B)<\delta}
 \int_BG(x,y)^{\,p}\mu(\d y)
&\leq
 \sup_{x\in E}
 \int_{B_r(x)}G(x,y)^{\,p}\mu(\d y)
+
 \sup_{B\subset B_R(o),\mu(B)<\delta}
 r^{p(\beta-\nu)}\mu(B)\\
&\leq
 \frac{\eps}{4}
+
 \frac{\eps}{4}
=
 \frac{\eps}{2}.
\end{align*}
Thus we have $\mu\in  S_{C\!K_{\infty}}^{\,p}({\bf X})$.
\end{apf}

\begin{apf}{Theorem~\ref{thm:Equals(1)} (1)}
We already know
$
 S_{C\!K_{\infty}}^{\,p}({\bf X}^{(1)})
\subset
 S_K^{\,p}({\bf X})
$.
Suppose that $\mu\in S_K^{\,p}({\bf X})$.
Throughout the proof, we fix
$\varepsilon > 0$, $\alpha > 0$ and $t\in ]0, t_0[$.

\vspace{0.5eM}

(Case I)
$\mu(E)<+\infty$ and $\nu \geq \beta$:
In this case,
by the upper bound of {\bf (A)}\ref{asmp:Phi},
we can see that for $d(z, y)\geq r$,
\begin{align*}
 \int_0 ^t
  p_s(z, y)
 \d s
 \leq
 \int_0 ^t
  \frac{1}{s^{\nu/\beta}}
  \Phi_2\left(\frac{r}{s^{1/\beta}}\right)
 \d s
\leq
 \beta r^{\beta-\nu}
 \int_{0} ^{\infty}
  u^{\nu-\beta-1}\Phi_2(u)
 \d u
=:
 M(r)
<
 +\infty
\end{align*}
and then, we have for any $a>0$,
\begin{align*}
 \int_{d(x, y)>r}
 \int_{d(z, y)\geq r}
  \left(
   \int_0 ^t
    p_s(z, y)
   \d s
  \right)^p
  p_{a}(x, z)
 \m (\d z)
 \mu(\d y)
\leq
 M(r)^p
 \mu(E)
.
\end{align*}
We also have
\begin{align*}
 \int_{d(x, y)>r}
 \int_{d(z, y)<r}
  \left(
   \int_0 ^t
    p_s(z, y)
   \d s
  \right)^p
  p_{a}(x, z)
 \m (\d z)
 \mu(\d y)
\leq
 \sup_{z\in E}
 \int_{B_r(z)}
  \left(
   \int_0 ^t
    p_s(z, y)
   \d s
  \right)^p
 \mu(\d y)
\end{align*}
and hence
\begin{align*}
 \int_{d(x, y)>r}
  \left(
   \int_{a} ^{a+t}
    p_s(z, y)
   \d s
  \right)^p
 \mu(\d y)
\leq
 M(r)^p
 \mu(E)
+
 \sup_{z\in E}
 \int_{B_r(z)}
  \left(
   \int_0 ^t
    p_s(z, y)
   \d s
  \right)^p
 \mu(\d y)
,
\end{align*}
which concludes that
\begin{align}
\notag
 \left\{
  \int_{d(x, y)>r}
   R_\alpha(x, y)^p
  \mu(\d y)
 \right\}^{\frac{1}{p}}
\leq&
 \sum_{n=0} ^\infty
 e^{-\alpha nt}
 \left\{
 \int_{d(x,y)>r}
  \left(
   \int_{nt} ^{(n+1)t}
    p_s(x, y)
   \d s
  \right)^p
 \mu (\d y)
 \right\}^{\frac{1}{p}}
\\
\label{eq:GreenEstimate_CaseI}
\leq&
 \frac{1}{1-e^{-\alpha t}}
 \left\{
  M(r)^p
  \mu(E)
 +
  \sup_{z\in E}
  \int_{B_r(z)}
   \left(
    \int_0 ^t
     p_s(z, y)
    \d s
   \right)^p
  \mu(\d y)
 \right\}^{\frac{1}{p}}
.
\end{align}
By \cite[Theorem~4.1]{KwM:pKato},
$\mu\in S_K^{\,p}({\bf X})$ implies
\begin{align*}
 \lim_{r\to0}
 \sup_{x\in E}
 \int_{B_r(x)}R_{\alpha}(x,y)^{\,p}\mu(\d y)
=
 0
\quad
\text{and}
\quad
 \lim_{r\to0}
 \sup_{x\in E}
 \int_{B_r(x)}
  \left(
   \int_0 ^t
    p_s(x, y)
   \d s
  \right)^p
 \mu(\d y)
=
 0
.
\end{align*}
We choose small $r>0$ so that
$
 \sup_{x\in E}
 \int_{B_r(x)}R_{\alpha}(x,y)^{\,p}\mu(\d y)
<
 \eps
$
and
$
 \sup_{x\in E}
 \int_{B_r(x)}
  \left(
   \int_0 ^t
    p_s(x, y)
   \d s
  \right)^p
 \mu(\d y)
<
 \eps
$,
a compact set $K$ satisfying
$
 M(r)^p
 \mu(K^c)
<
 \eps
$
and $\delta>0$ so that
$
 M(r)^p
 \delta
<
 \eps
$.
Then applying \eqref{eq:GreenEstimate_CaseI}
by replacing $\mu$ with $\1_{K^c\cup B} \mu$,
we have for any Borel set $B\subset K$ with $\mu(B)<\delta$
\begin{align*}
 \sup_{x\in E}&
 \int_{K^c\cup B}R_{\alpha}(x,y)^{\,p}\mu(\d y)
\\
\leq&
 \sup_{x\in E}
 \int_{B_r(x)}R_{\alpha}(x,y)^{\,p}\mu(\d y)
+
 \sup_{x\in E}
 \int_{B_r(x)^c}
  R_{\alpha}(x,y)^{\,p}
  \1_{K^c\cup B}(y)
 \mu(\d y)
\\
\leq&
 \left\{
  1
 +
  3
  \left(
   \frac{1}{1-e^{-\alpha t}}
  \right)^p
 \right\}
 \eps
.
\end{align*}
So we can obtain
$\mu\in S_{C\!K_{\infty}}^{\,p}({\bf X}^{(1)})$ for this case.

\vspace{0.5eM}

(Case I\hspace{-0.1eM}I)
$\mu(E)<+\infty$ and $\nu < \beta$:
In this case, we see that
\begin{align*}
 \int_0 ^t
  p_s(z, y)
 \d s
\leq
 \Phi_2(0)
 \int_0 ^t
  \frac{1}{s^{\nu/\beta}}
 \d s
<
 +\infty
.
\end{align*}
By a similar calculation as to obtain
\eqref{eq:GreenEstimate_CaseI},
we have
\begin{align*}
 \sup_{x\in E}
 \int_{K^c\cup B}R_{\alpha}(x,y)^{\,p}\mu(\d y)
&\leq
 \left(
  \frac{1}{1-e^{-\alpha t}}
 \right)^p
 \left(
  \Phi_2(0)
  \int_0 ^t
   \frac{1}{s^{\nu/\beta}}
  \d s
 \right)^{\,p}
 \mu(K^c\cup B).
\end{align*}
So we can obtain
$\mu\in S_{C\!K_{\infty}}^{\,p}({\bf X}^{(1)})$
for this case.

\vspace{0.5eM}

(Case I\hspace{-0.1eM}I\hspace{-0.1eM}I)
$\mu \in S_{K_{\infty}}^1({\bf X}^{(1)})$ with $p>1$
and $\nu \geq \beta$:
In this case, we have for any $a>0$,
\begin{align*}
 \int_{d(x, y)>r}
 \int_{d(z, y)\geq r}&
  \left(
   \int_0 ^t
    p_s(z, y)
   \d s
  \right)^p
  p_{a}(x, z)
 \m (\d z)
 \mu(\d y)
\\
\leq&
 M(r)^{p-1}
 \sup_{z\in E}
 \int_E
  \left(
   \int_0 ^t
    p_s(z, y)
   \d s
  \right)
 \mu (\d y)
\\
\leq&
 M(r)^{p-1}
 e^{\alpha t}
 \sup_{z\in E}
 \int_E
  R_\alpha (z, y)
 \mu (\d y)
.
\end{align*}
By a similar calculation
to obtain \eqref{eq:GreenEstimate_CaseI},
we have
\begin{align}
\notag&\hspace{-0.6cm}
 \left\{
  \int_{d(x, y)>r}
   R_\alpha(x, y)^p
  \mu(\d y)
 \right\}^{\frac{1}{p}}
\\
\label{eq:GreenEstimate_CaseIII}
\leq&
 \frac{1}{1-e^{-\alpha t}}
 \left\{
  M(r)^{p-1}
  e^{\alpha t}
  \sup_{z\in E}
  \int_E
   R_\alpha (z, y)
  \mu (\d y)
 +
  \sup_{z\in E}
  \int_{B_r(z)}
   \left(
    \int_0 ^t
     p_s(z, y)
    \d s
   \right)^p
  \mu(\d y)
 \right\}^{\frac{1}{p}}
.
\end{align}
Now, choose small $r>0$ so that
$
 \sup_{x\in E}
 \int_{B_r(x)}R_{\alpha}(x,y)^{\,p}\mu(\d y)
<
 \eps
$
and
$
 \sup_{x\in E}
 \int_{B_r(x)}
  \left(
   \int_0 ^t
    p_s(x, y)
   \d s
  \right)^p
 \mu(\d y)
<
 \eps
$,
a compact set $K$ satisfying
$
 \allowbreak
 M(r)^{p-1}
 e^{\alpha t}
 \sup_{x\in E}
 \int_{K^c}R_{\alpha}(x,y)\mu(\d y)
<
 \eps
$
and $\delta>0$ so that
$
 M(r)^{p}
 \delta
<
 \eps
$.
Apply
\eqref{eq:GreenEstimate_CaseI}
by replacing $\mu$ with $\1_{B} \mu$,
and apply
\eqref{eq:GreenEstimate_CaseIII}
by replacing $\mu$ with $\1_{K^c} \mu$.
Then
we have for any Borel set $B\subset K$ with $\mu(B)<\delta$
\begin{align*}
&\hspace{-0.4cm}
 \sup_{x\in E}
 \int_{K^c\cup B}R_{\alpha}(x,y)^{\,p}\mu(\d y)
\\
\leq&
 \sup_{x\in E}
 \int_{B_r(x)}R_{\alpha}(x,y)^{\,p}\mu(\d y)
+
 \sup_{x\in E}
 \int_{B_r(x)^c}
  R_{\alpha}(x,y)^{\,p}
  \1_{B}(y)
 \mu(\d y)
+
 \sup_{x\in E}
 \int_{B_r(x)^c}
  R_{\alpha}(x,y)^{\,p}
  \1_{K^c}(y)
 \mu(\d y)
\\
\leq&
 \left\{
  1
 +
  4
  \left(
   \frac{1}{1-e^{-\alpha t}}
  \right)^p
 \right\}
 \eps.
\end{align*}
So we can obtain
$\mu\in S_{C\!K_{\infty}}^{\,p}({\bf X}^{(1)})$ for this case.

\vspace{0.5eM}

(Case I\hspace{-0.1eM}V)
$\mu \in S_{K_{\infty}}^1({\bf X}^{(1)})$ with $p>1$
and $\nu < \beta$:
In this case, we see
\begin{align*}
 \int_0 ^t
  p_s(z, y)
 \d s
\leq
 \Phi_2(0)
 \int_0 ^t
  \frac{1}{s^{\nu/\beta}}
 \d s
<
 +\infty
.
\end{align*}
By a similar calculation as to obtain
\eqref{eq:GreenEstimate_CaseIII},
we have
\begin{align*}
 \sup_{x\in E}
 \int_{K^c}R_{\alpha}(x,y)^{\,p}\mu(\d y)
&\leq
 \left(
  \frac{1}{1-e^{-\alpha t}}
 \right)^p
 \left(
  \Phi_2(0)
  \int_0 ^t
   \frac{1}{s^{\nu/\beta}}
  \d s
 \right)^{p-1}
 e^{\alpha t}
 \sup_{x\in E}
 \int_{K^c}R_{\alpha}(x,y)\mu(\d y)
.
\end{align*}
So we can obtain
$\mu\in S_{C\!K_{\infty}}^{\,p}({\bf X}^{(1)})$
for this case
by the same manner as
(Case I\hspace{-0.1eM}I\hspace{-0.1eM}I).
\end{apf}

\begin{apf}{Theorem~\ref{thm:Equals(1)} (2)}
Since
$
 S_{C\!K_{\infty}}^{\,p}({\bf X}^{(1)})
\subset
 S_{K_{\infty}}^{\,p}({\bf X}^{(1)})
$
 by Remark~\ref{rem:GreenTightInclusion}, it suffices to show
$
 S_{K_{\infty}}^{\,p}({\bf X}^{(1)})
\subset
 S_{C\!K_{\infty}}^{\,p}({\bf X}^{(1)})
$.
Take
$
 \mu
\in
 S_{K_{\infty}}^{\,p}({\bf X}^{(1)})
\subset
 S_K^{\,p}({\bf X})
$.
Then, for any $\eps>0$ there exists a compact set $K$ such that
$
 \sup_{x\in E}
 \int_{K^c}
  R_1(x,y)^p
 \mu(\d y)
<
 \eps/2
$.
For given $\delta>0$,
we take any Borel subset $B$ of $K$ with $\mu(B)<\delta$.
First we assume $\nu\geq\beta$.
By applying
\cite[Theorem~4.1]{KwM:pKato} with $\mu$,
for any $\eps>0$, there exists $r>0$ such that
\begin{align*}
 \sup_{z\in E}
 \int_{B_r(z)}
  R_1(z,y)^p
 \mu(\d y)
<
 \eps
\quad \text{ and }\quad
 \sup_{z\in E}
 \int_{B_r(z)}
  \left(
   \int_0^tp_s(z,y)\d s
  \right)^{p}
 \mu(\d y)
<
 \frac{\eps}{2},
\end{align*}
where $t\in]0,t_0[$ is a fixed small time.
Take a small $r>0$ so that $M(r)^p\delta<\eps/2$.
Applying \eqref{eq:GreenEstimate_CaseI} with $\alpha=1$
by replacing $\mu$ with $\1_B\mu$,
\begin{align*}
 \left(
  \int_{B\cap B_r(x)^c}
   R_1(x,y)^p
  \mu(\d y)
 \right)^{\frac{1}{p}}
\leq
 \frac{1}{1-e^{-t}}
 \left\{
  M(r)^p\mu(B)
 +
  \sup_{z\in E}
  \int_{B_r(z)}
   \left(
    \int_0^tp_s(z,y)\d s
   \right)^{\,p}
  \mu(\d y)
 \right\}^{\frac{1}{p}}.
\end{align*}
Then we have
\begin{align*}
 \int_B
  R_1(x,y)^{\,p}
 \mu(\d y)
&\leq
 \int_{B\cap B_r(x)}
  R_1(x,y)^{\,p}
 \mu(\d y)
+
 \int_{B\cap B_r(x)^c}
  R_1(x,y)^{\,p}
 \mu(\d y)
\\
&\leq
 \eps
+
 \left(
  \frac{1}{1-e^{-t}}
 \right)^{p}
 \left(
  \frac{\eps}{2}
 +
  \frac{\eps}{2}
 \right)
=
 \left(
  1
 +
  \left(\frac{1}{1-e^{-t}}\right)^{p\,}
 \right)
 \eps.
\end{align*}
This concludes $\mu\in S_{C\!K_{\infty}}^{\,p}({\bf X}^{(1)})$.

Next we suppose $\nu<\beta$.
As noted in the proof of Theorem~\ref{thm:Equals(1)} (1)
(Case I\hspace{-0.1em}I),
\begin{align*}
 \int_0^t
  p_s(z,y)
 \mu(\d y)
\leq
 \Phi_2(0)
 \int_0^t
  \frac{1}{s^{\nu/\beta}}
 \d s
<
 +\infty.
\end{align*}
By a similar calculation as to obtain \eqref{eq:GreenEstimate_CaseI}
with $\alpha=1$, we have
\begin{align*}
 \sup_{x\in E}
 \int_B
  R_1(x,y)^{\,p}
 \mu(\d y)
\leq
 \left(
  \frac{1}{1-e^{-t}}
 \right)^{\,p}
 \left(
  \Phi_2(0)
  \int_0^t
   \frac{1}{s^{\nu/\beta}}
  \d s
 \right)^{\,p}
 \mu(B).
\end{align*}
So we obtain $\mu\in S_{C\!K_{\infty}}^{\,p}({\bf X}^{(1)})$
for this case.
\end{apf}

\section{Examples}\label{sec:Examples}

\begin{ex}[{Brownian motions on \boldmath$\R^d$}]\label{ex:Brownian}
{\rm
Let ${\bf X}=(\Omega, B_t, {\PP}_x)_{x\in\R^d}$
be a $d$-dimensional Brownian motion on $\R^d$.
Consider $p\in[1,+\infty[$.
We say that $\mu\in K_{d}^{\,p}$ (or $\mu\in K_{d,2}^{\,p}$)
if and only if
\begin{align*}
 \lim_{r\to0}
 \sup_{x\in\R^d}
 \int_{|x-y|<r}\frac{\mu(\d y)}{|x-y|^{(d-2)p}}
=
 0
\quad &\text{ for }\quad d\geq3,\\
 \lim_{r\to0}
 \sup_{x\in\R^d}
 \int_{|x-y|<r}(\log|x-y|^{-1})^p\mu(\d y)
=
 0
\quad &\text{ for }\quad d=2, \\
 \sup_{x\in\R^d}
 \int_{|x-y|\leq1}\mu(\d y)
<
 +\infty
\quad &\text{ for }\quad d=1.
\end{align*}
We write $K_d$ instead of $K_d^1$ for $p=1$.
Then we have $K_d^{\,p}=S_K^{\,p}({\bf X})$ by
\cite[Example~2.4]{TM:pKato}.
The $d$-dimensional Lebesgue measure $\m$ belongs to
$K_d^{\,p}=S_K^{\,p}({\bf X})$ if and only if
$p\in [1, d/(d-2)_+[$
by \cite[Theorem~3.2 or Corollary~4.4]{KwM:pKato},
where
$d/(d-2)_+ := d/(d-2)$ if $d\geq 3$,
$d/(d-2)_+ := +\infty$ if $d=1, 2$.
For any non-negative bounded $g\in L^1(\R^d)$
the finite measure $g\m$ also belongs to
$S_{C\!K_{\infty}}^{\,p}({\bf X}^{(1)})$
(to $S_{C\!K_{\infty}}^{\,p}({\bf X})$
if ${\bf X}$ is transient)
by Theorem~\ref{thm:Equals}
under $p\in [1, d/(d-2)_+[$.
The surface measure $\sigma_R$
on the $R$-sphere $\partial B_R(0)$ satisfies that
$\sigma_R(B_r(x))\leq C_2r^{d-1}$
for any $x\in\R^d$ and $r>0$ with some $C_2>0$,
and
$\sigma_R(B_r(x))\geq C_1r^{d-1}$
for any $x\in \partial B_R(0)$ and $r\in ]0,r_0[$
with some $C_1, r>0$.
Then
we can conclude that by Theorem~\ref{thm:Equals} and
\cite[Theorem~3.2]{KwM:pKato},
$\sigma_R\in S_{C\!K_{\infty}}^{\,p}({\bf X}^{(1)})$ holds
if and only if $p\in [1, (d-1)/(d-2)_+[$,
where
$(d-1)/(d-2)_+ := (d-1)/(d-2)$ if $d\geq 3$,
$(d-1)/(d-2)_+ := +\infty$ if $d=1,2$.
Moreover,
$\sigma_R\in S_{C\!K_{\infty}}^{\,p}({\bf X})$ holds
if and only if $p\in [1, (d-1)/(d-2)[$
provided $d\geq3$.

We consider a non-empty connected open set $D$ of $\R^d$.
The boundary point $z\in \partial D$ is said to be
\emph{regular} if ${\PP}_z(\tau_D=0)=1$.
Denote by $(\partial D)_r$ the set of regular points
in boundary.
$D$ is said to be \emph{regular} if
$(\partial D)_r=\partial D$.
Let $D$ be a connected open regular set of $\R^d$.
The absorbing Brownian motion
${\bf X}_D=(\Omega, X_t^D,{\PP}_x)$
(or part process of ${\bf X}$ on $D$) is defined
as the process killed upon leaving $D$.
Then
${\bf X}_D$ is an irreducible doubly Feller diffusion process
on $D$ (see \cite{Chungdoubl}).
If further $D^c$ is non-polar, (in particular $\m(D^c)>0$),
then ${\bf X}_D$ is transient
in view of \cite[Theorem~4.7.1 and Exercise~4.7.1]{FOT}.
Let $R^D(x,y)$ be the Green function with respect to
${\bf X}_D$.
$D$ is said to be \emph{Green-bounded} if
$
 \sup_{x\in D}
 \int_D R^D(x,y){\rm d}y
=
 \sup_{x\in D}
 {\EE}_x[\tau_D]
<
 +\infty
$,
equivalently $\m\in S_{D_0}^1({\bf X}_D)$,
where $\m$ is the $d$-dimensional Lebesgue measure on $D$.
If $d=1$ and $D$ is not bounded,
$R^D\nu\in C_{\infty}(D)$ fails even for $\nu(D)<+\infty$
(see \cite[Example~1]{KimKuwae:GeneralAnal}).
Since
$
 \sup_{x\in D}
 {\EE}_x[\tau_D]
\leq
 \frac{d+2}{2\pi d}
 \left(\frac{d+2}{2}\right)^{2/d}
 \m(D)^{2/d}
$
(see \cite[Theorem~1.17]{CZ:BS}), $\m(D)<+\infty$,
in particular the boundedness of $D$,
implies the Green-boundedness of $D$.
For a (positive) Radon measure $\nu$ on $\R^d$, if
\begin{align}
 \left\{
 \begin{array}{cl}
  d=1, & D\text{ is bounded and }\nu(D)<+\infty\text{ or, } \\
  d=2, & D\text{ is Green-bounded and } \nu\in S_K^{\,p}({\bf X})\text{ with }\nu(D)<+\infty
\text{ or, } \\
  d\geq3, & \nu\in S_K^{\,p}({\bf X})\text{ with }\nu(D)<+\infty,
 \end{array}
 \right.
\label{eq:KatoGreen}
\end{align}
then we can prove
$
 (R^D)^{\,p}\nu
:=
 \int_D R^D(\cdot,y)^{\,p}\nu(\d y)
\in
 C_b(D)
$,
and it belongs to $C_{\infty}(D)$
provided $D$ is a regular domain.
Indeed, by \cite[Theorem~2.6(ii)]{CZ:BS},
\begin{align}
 R^D(x,y)
\leq
 \left\{
 \begin{array}{cc}
  \displaystyle{\frac{1}{\pi}\log^+|x-y|^{-1}+C} & d=2, \\
  \displaystyle{C|x-y|^{-(d-2)}}& d\geq3,
 \end{array}
 \right.
\label{eq:GreenKernelEst}
\end{align}
where $C$ is a positive constant depending on
$\|R^D\m\|_{\infty}$ for $d=2$
(here we use the Green-boundedness of $D$)
and on $d$ for $d\geq3$.
Moreover, there exists a positive sequence $\alpha_n\to0$
such that
$R^D(x,y)=R^D(x,y)\land n$ if $|x-y|\geq\alpha_n$.
Then we can calculate
\begin{align}
 \sup_{x\in D}
 \left|
  \int_DR^D(x,y)^{\,p}\nu(\d y)
 -
  \int_D(R^D(x,y)\land n)^{\,p}\nu(\d y)
 \right|
\leq
 \sup_{x\in D}
 \int_{\{|x-y|<\alpha_n\}}R^D(x,y)^{\,p}\nu(\d y).
\label{eq:GreenPB}
\end{align}
The right-hand side of \eqref{eq:GreenPB}
uniformly converges to $0$ as $n\to\infty$,
because of the estimate
\eqref{eq:GreenKernelEst},
$\nu\in S_{K}^{\,p}({\bf X})=K_d^{\,p}$ and $\nu(D)<+\infty$
for the case $d=2$.
Under the conditions in \eqref{eq:KatoGreen}, the function
\begin{align*}
 x
\mapsto
 \int_D(R^D(x,y)\land n)^{\,p}\nu(\d y)
\end{align*}
belongs to $C_b(D)$ and in $C_{\infty}(D)$
provided $D$ is a regular domain,
because of the extended continuity of
$(x,y)\mapsto R^D(x,y)$ (see \cite[Theorem~2.6(iii)]{CZ:BS}).
The uniform convergence
\begin{align*}
 \lim_{n\to\infty}
 \lim_{x\in D}
 \left|
  \int_D R^D(x,y)^{\,p}\nu(\d y)
 -
  \int_D(R^D(x,y)\land n)^{\,p}\nu(\d y)
 \right|
=
 0
\end{align*}
noted above implies the assertion for the case $d\geq2$.
The proof of $(R^D)^p\nu\in C_{\infty}(D)$ for $d=1$ is
clear from the expression
\begin{align}
 \int_a^b R^D(x,y)^{\,p}\nu(\d y)
=
 \left(\frac{2(x-a)}{b-a} \right)^p
 \int_x^b(b-y)^p\nu(\d y)
+
 \left(\frac{2(b-x)}{b-a} \right)^p
 \int_a^x(y-a)^p\nu(\d y)
\label{eq:ExpressionOneDim}
\end{align}
for $D=]\,a,\,b\,[$.
Hence, if $\nu$ satisfies \eqref{eq:KatoGreen} and
$D$ is a regular domain, then one can obtain
$\1_D \nu\in C_\infty(D)$.
For any compact set $K$ of $D$,
the $0$-order version of Frostman's maximum principle
\eqref{eq:Frostman} gives that
$
 \sup_{x\in D}
 (R^D)^{\,p}\1_{K^c}\nu (x)
=
 \sup_{x\in D\setminus K}
 (R^D)^{\,p}\1_{K^c}\nu (x)
\leq
 \sup_{x\in D\setminus K}
 (R^D)^{\,p}\nu (x)
$.
Therefore
\begin{align*}
 \1_D\nu
\in
 S_{K_{\infty}}^{\,p}({\bf X}_D)
=
 S_{C\!K_{\infty}}^{\,p}({\bf X}_D),
\end{align*}
where the equality follows from Theorem~\ref{thm:Coincidence},
because ${\bf X}_D$ possesses {\bf (RSF)}.
It is proved in \cite[Corollary~3.3]{KwM:pKato} that
$p(d-2)<d$ is equivalent to $\m\in S_K^{\,p}({\bf X})$.
From this, if
\begin{align}
 \left\{
 \begin{array}{cl}
  d=1, & D\text{ is bounded}\text{ or, }\\
  d\geq2, & p(d-2)<d\text{ with }\m(D)<+\infty,
 \end{array}
 \right.
\label{eq:KatoGreenUnderlying}
\end{align}
then
$
 \1_D\m
\in
 S_{K_{\infty}}^{\,p}({\bf X}_D)
=
 S_{C\!K_{\infty}}^{\,p}({\bf X}_D)
$
provided $D$ is a regular domain.
We relax \eqref{eq:KatoGreenUnderlying} in the following way:
\begin{align}
 \left\{
 \begin{array}{cl}
  d=1, & D\text{ is bounded}\text{ or, } \\
  d\geq2, & p(d-2)<d\text{ with }\lim_{x\in D,|x|\to \infty}\m(D\cap B_1(x))=0.
 \end{array}
 \right.
\label{eq:KatoGreenUnderlyingRelaxed}
\end{align}
Therefore, by Theorem~\ref{thm:compactEmbedding},
we have the following:
\begin{thm}\label{thm:CptEmbedEuclideanRelaxed}
Suppose that
\eqref{eq:KatoGreenUnderlyingRelaxed} is satisfied.
Then
$
 \1_D\m
\in
 S_{K_{\infty}}^{\,p}({\bf X}_D^{(1)})
=
 S_{C\!K_{\infty}}^{\,p}({\bf X}_D^{(1)})
$,
hence the embedding
\begin{align*}
 H_0^1(D)
\hookrightarrow
 L^{2p}(D;\m)
\end{align*}
is compact.
Moreover, if $D$ is Green-bounded, then
$
 \1_D\m
\in
 S_{K_{\infty}}^{\,p}({\bf X}_D)
=
 S_{C\!K_{\infty}}^{\,p}({\bf X}_D)
$,
hence the embedding
\begin{align*}
 H_0^1(D)_e
\hookrightarrow
 L^{2p}(D;\m)
\end{align*}
is compact.
\end{thm}
\begin{pf}
The latter assertion follows from
Proposition~\ref{prop:GreenBddEqui}.
So it suffices to prove the former assertion.
Since ${\bf X}_D$ possesses {\bf (RSF)},
we know
$
 \1_D\m
\in
 S_{K_{\infty}}^1({\bf X}_D^{(1)})
=
 S_{C\!K_{\infty}}^1({\bf X}_D^{(1)})
$
under \eqref{eq:KatoGreenUnderlyingRelaxed}
by \cite[Lemma~3.3]{TakTawTsuchi:Compact}.
So there exists an increasing sequence $\{K_{\ell}\}$
of compact subsets of $D$ such that
\begin{align}
 \lim_{\ell\to\infty}
 \sup_{x\in D}
 \int_{D\setminus K_{\ell}}R_1^D(x,y)\m(\d y)
=
 0.
\label{eq:GreentightOne}
\end{align}
It is easy to see that there exist $C=C(d)>0$ such that
\begin{align}
 R_1^D(x,y)
\leq
 R_1(x,y)
\leq
 \left\{\begin{array}{cc}
  e^{-\sqrt{2}|x-y|}/\sqrt{2} & d=1, \\
  \displaystyle{
   \frac{1}{\pi}\cdot\frac{1}{|x-y|^2}+\frac{1}{2\pi}
  }
  & d=2, \\
  \displaystyle{\frac{C}{|x-y|^{d-2}}} & d\geq3,
 \end{array}
 \right.
\qquad\text{ for \ \ \ all }\qquad x,y\in D.
\label{eq:GreenEstimates}
\end{align}
Indeed, for $d=2$,
\begin{align*}
 R_1(x,y)
&\leq
 \int_0^{\frac{|x-y|^2}{2}}
  e^{-t}\frac{1}{2\pi t}
  e^{-\frac{|x-y|^2}{2t}}
 \d t
+
 \int_{\frac{|x-y|^2}{2}}^{\infty}
  e^{-t}\frac{1}{2\pi t}
  e^{-\frac{|x-y|^2}{2t}}
 \d t\\
&\leq
 \int_1^{\infty}
  e^{-s}\frac{1}{2\pi s}
  e^{-\frac{|x-y|^2}{2s}}
 \d s
+
 \frac{1}{2\pi}\cdot\frac{2}{|x-y|^2}
 \int_{\frac{|x-y|^2}{2}}^{\infty}e^{-t}\d t\\
&\leq
 \frac{1}{2\pi}
+
 \frac{1}{\pi}\cdot\frac{1}{|x-y|^2}.
\end{align*}
Then there exists a decreasing sequence $\{\alpha_n\}$
converging to $0$ such that
$R_1^D(x,y)=R_1^D(x,y)\land n$ for $|x-y|\geq\alpha_n$.
Indeed,
we can choose
$\alpha_n = (C_3/n)^{\frac{1}{d-2}}$ for $d\geq3$,
$\alpha_n = \sqrt{\frac{2}{2\pi n-1}}$ for $d=2$
with $n>1/(2\pi)$, and $\alpha_n$ is arbitrary for $d=1$.
On the other hand,
\begin{align}
 \hspace{-1cm}
 \left|
  \sup_{x\in D}
  \int_{D\setminus K_{\ell}}
   R_1^D(x,y)^{\,p}
  \m(\d y)
 -
  \sup_{x\in D}
  \int_{D\setminus K_{\ell}}
   R_1^D(x,y)
   \left(R_1^D(x,y)\land n \right)^{p-1}
  \m(\d y)
 \right|
&\leq
 \sup_{x\in \R^d}
 \int_{|x-y|<\alpha_n}
  R_1(x,y)^{\,p}
 \m(\d y).
\label{eq:LpKatoUnif}
\end{align}
Since $\m\in S_{K}^{\,p}({\bf X})=K_{d,2}^{\,p}$,
the right-hand side of \eqref{eq:LpKatoUnif} converges to $0$
by \cite[Theorem~4.1]{KwM:pKato}.
Combining this with \eqref{eq:GreentightOne}, we see
\begin{align*}
 \lim_{\ell\to\infty}
 \sup_{x\in D}
 \int_{D\setminus K_{\ell}}
  R_1^D(x,y)^{\,p}
 \m(\d y)
=
 0,
\end{align*}
that is, $\1_D\m\in S_{K_{\infty}}^{\,p}({\bf X}_D^{(1)})$.
For $d=1$ with bounded $D$, it is easy to see
the same assertion by
$
 \1_D\m
\in
 S_{K_{\infty}}^{\,p}({\bf X}_D)
=
 S_{C\!K_{\infty}}^{\,p}({\bf X}_D)
$
derived from the expression \eqref{eq:ExpressionOneDim}
by replacing $\nu$ with $\m$.
\end{pf}

Next we set
\begin{align}
 \mathscr{B}_0
:=
 \left\{
  B\in\mathscr{B}(\R^d)
 \;\left|\;
  \lim_{|x|\to\infty}\m(B\cap B_1(x))=0
 \right.
 \right\}.
\label{eq:B0}
\end{align}
As noted in \cite[(4.1)]{TakTawTsuchi:Compact},
it holds that
\begin{align}
 \lim_{|x|\to\infty}
 \m(B\cap B_R(x))
=
 0
\quad \text{ for \ \ \ any }\quad R>0.
\label{eq:(4.1)}
\end{align}
As proved in \cite[Theorem~4.1]{TakTawTsuchi:Compact},
a domain $B$ belongs to $\mathscr{B}_0$
if and only if
$
 \m^B
:=
 \1_B\m
\in
 S_{K_{\infty}}^1({\bf X}^{(1)})
=
 S_{C\!K_{\infty}}^1({\bf X}^{(1)})
$
and if and only if
$H^1(\mathbb{R}^d)$ is compactly embedded into $L^2(D)$
provided $d\geq3$.
In the following,
we will give a $p$-extension of the fact.

\begin{lem}\label{lem:GrennTightBM}
For general $d\geq1$,
$B\in\mathscr{B}_0$ implies
$
 \m^B
\in
 S_{K_{\infty}}^1({\bf X}^{(1)})
=
 S_{C\!K_{\infty}}^1({\bf X}^{(1)})
$.
\end{lem}
\begin{pf}
By the ultra-contractivity
$
 \|P_t\|_{L^1\rightarrow L^\infty}
\leq
 C
 t^{-\frac{d}{2}}
$
of ${\bf X}$ and the estimate
\begin{align*}
 \int_0^tP_s\1_{B\setminus B_R(x)}(x)\d s
\leq&
 {\EE}_x\left[\int_0^t\1_{B_R(x)^c}(X_s)\d s \right]\\
\leq&
 {\EE}_x[(t-\tau_{B_R(x)})_+]
=
 {\EE}_0[(t-\tau_{B_R(0)})_+]\\
\leq&
 t{\PP}_0(\tau_{B_R(0)}\leq t),
\end{align*}
we have
\begin{align*}
 \int_0^t P_s\1_B(x)\d s
&\leq
 \int_0^t P_s\1_{B\cap B_R(x)}(x)\d s
+
 \int_0^tP_s\1_{B\setminus B_R(x)}(x)\d s\\
&\leq
 C
 \int_0^t s^{-\frac{d}{2q}}\d s
\cdot
 \m(B\cap B_R(x))^{\frac{1}{q}}
+
 t{\PP}_0(\tau_{B_R(0)}\leq t),
\end{align*}
where $q>d/2$.
So
\begin{align*}
 \varlimsup_{|x|\to\infty}
 \int_0^{t}P_{s+nt}\1_B(x)\d s
&\leq
 (t+nt){\PP}_0(\tau_{B_R(0)}\leq t+nt)
\longrightarrow
 0
\quad\text{ as }\quad R\to\infty.
\end{align*}
Here we use the quasi-left continuity of ${\bf X}$.
From this,
\begin{align*}
 R_1\1_B(x)
&=
 \sum_{n=0}^{\infty}
 e^{-nt}
 \int_0^t e^{-s}P_{s+nt}\1_B(x)\d s\\
&\leq
 \sum_{n=0}^{\infty}
 e^{-nt}
 \int_0^t e^{-s}P_{s+nt}\1_B(x)\d s
\longrightarrow
 0
\quad\text{ as }\quad |x|\to\infty,
\end{align*}
which implies $R_1\1_B\in C_{\infty}(\R)$, hence
$
 \m^B
\in
 S_{K_{\infty}}^1({\bf X})
=
 S_{C\!K_{\infty}}^1({\bf X})
$.
\end{pf}
Now we claim the following:
\begin{prop}\label{prop:compactEmbedEuc}
Suppose $p(d-2)<d$ and $B\in\mathscr{B}_0$.
Then
$
 \m^B
\in
 S_{K_{\infty}}^{\,p}({\bf X}^{(1)})
=
 S_{C\!K_{\infty}}^{\,p}({\bf X}^{(1)})
$.
In particular,
$H^1(\R^d)$ is compactly embedded into $L^{2p}(B)$.
\end{prop}
\begin{pf}
It suffices to prove
$\m^B\in S_{K_{\infty}}^{\,p}({\bf X}^{(1)})$ by
Theorems~\ref{thm:compactEmbedding} and \ref{thm:Coincidence}.
By \eqref{eq:GreenEstimates}, we have
$R_1(x,y)=R_1(x,y)\land n$ for $|x-y|\geq\alpha_n$,
where $\alpha_n$ is the constant appeared as above.
Then one can deduce the following estimate:
\begin{align}
 \left|
  \sup_{x\in\R^d}
  \int_{K_{\ell}^c}\hspace{-0.2cm}R_1(x,y)^{\,p}\m^B(\d y)
 -
  \sup_{x\in\R^d}
  \int_{K_{\ell}^c}\hspace{-0.2cm}
   R_1(x,y)(R_1(x,y)\land n)^{p-1}
  \m^B(\d y)
 \right|
\leq
 \sup_{x\in\R^d}
 \int_{|x-y|<\alpha_n}\hspace{-0.5cm}
  R_1(x,y)^{\,p}
 \m(\d y).
\label{eq:GreenEstimateBM}
\end{align}
The right-hand side of \eqref{eq:GreenEstimateBM}
converges to $0$ as $n\to\infty$
from $\m\in S_K^{\,p}({\bf X})$ under $p(d-2)<d$
by applying \cite[Theorem~3.1]{KwM:pKato} provided $d\geq2$.
When $d=1$, the right-hand side of \eqref{eq:GreenEstimateBM}
is estimated above by $(1/\sqrt{2})^p\m(B_{\alpha_n}(0))$,
which goes to $0$ as $n\to\infty$.
Thus we can obtain the assertion.
\end{pf}

\begin{thm}\label{thm:EdmundEvance}
Let $D$ be a domain of $\R^d$.
Suppose $p(d-2)<d$.
Then the following statements are equivalent:
\begin{enumerate}
\item
$D\in\mathscr{B}_0$,\quad
\item
$
 \1_D\m
\in
 S_{K_{\infty}}^{\,p}({\bf X}^{(1)})
=
 S_{C\!K_{\infty}}^{\,p}({\bf X}^{(1)})
$,\quad
\item
$H^1(\R^d)$ is compactly embedded into $L^{2p}(D)$.
\end{enumerate}
\end{thm}
\begin{pf}
The condition $p(d-2)<d$ is only used to establish
the continuity of the embedding
$H^1(\R^d)\hookrightarrow L^{2p}(D)$.
We have already proved
(1)$\Longrightarrow$(2)$\Longrightarrow$(3)
in Proposition~\ref{prop:compactEmbedEuc}.
The proof of the implication
(3)$\Longrightarrow$(1) is similar to
\cite[Chapter X, Lemma~6.11]{EdEv:Spectral}.
\end{pf}

In the end of this example,
we give the compactness of the Schr\"odinger semigroups,
which is a $p$-version of
\cite[Theorem~5.2]{TakTawTsuchi:Compact} under $\alpha=2$.
\begin{thm}\label{thm:SchoedingerCopact}
Let $V$ be a positive Borel function on $\R^d$
satisfying $V\m\in S_{L\!K}^1({\bf X})$.
Suppose that $\{V\leq M\}\in\mathscr{B}_0$ for any $M>0$ and
$p(d-2)<d$.
Then the Schr\"odinger semigroup $P_t^{-V}$ defined by
\begin{align}
 P_t^{-V}f(x)
=
 {\EE}_x\left[e^{-\int_0^tV(X_s)\d s}f(X_t)\right],
\quad
 f\in L^2(\R^d)\cap \mathscr{B}(\R^d),
\label{eq:SchroedingerSemigroup}
\end{align}
forms a compact operator from $L^2(\R^d)$ to $L^{2p}(\R^d)$.
\end{thm}
\begin{cor}\label{cor:SchoedingerCopact}
Let $V$ be a positive continuous function on $\R^d$
satisfying $V\m\in S_{L\!K}^1({\bf X})$.
Suppose that $p(d-2)<d$ and
\begin{align}
 \lim_{|x|\to\infty}
 V(x)
=
 +\infty.
\label{eq:exhaustion}
\end{align}
Then the Schr\"odinger semigroup $P_t^{-V}$
defined in \eqref{eq:SchroedingerSemigroup}
forms a compact operator from $L^2(\R^d)$ to $L^{2p}(\R^d)$.
\end{cor}
\begin{pf}
For any $M>0$, the sublevel set $\{V\leq M\}$ is a compact set
from \eqref{eq:exhaustion},
hence it belongs to $\mathscr{B}_0$ automatically.
\end{pf}

\begin{rem}\label{rem:Matsuura}
{\rm
Theorem~\ref{thm:SchoedingerCopact} and
Corollary~\ref{cor:SchoedingerCopact}
are not included in \cite[Example~4.4]{MatsuKohei:compact}.
Our conclusion of the compactness of $P_t^{-V}$ is different
from that in \cite[Theorems~2.2 and 2.4]{MatsuKohei:compact}.
}
\end{rem}

\begin{apf}{Theorem~\ref{thm:SchoedingerCopact}}
When $p=1$, the assertion is nothing but
\cite[Theorem~5.2]{TakTawTsuchi:Compact},
which was done by showing
$
 \m
\in
 S_{K_{\infty}}^1({\bf X}^{-V-1})
=
 S_{C\!K_{\infty}}^1({\bf X}^{-V-1})
$
(note that its proof is valid for general $d\geq 1$).
Here ${\bf X}^{-V-1}$ is the subprocess of ${\bf X}$
by $\exp\left(-\int_0^tV(X_s)\d s-t \right)$,
which possesses {\bf (RSF)}
because of $V\m\in S_{L\!K}^1({\bf X})$
(see \cite[Corollary~6.1]{KKK:Kyoto}).
Let $R_1^{-V}(x,y)$ be the Green kernel of ${\bf X}^{-V-1}$.
Then
\begin{align}
 \left|
  \sup_{x\in\R^d}
  \int_{K_{\ell}^c}\hspace{-0.2cm}
   R_1^{-V}(x,y)^{\,p}
  \m(\d y)
 -
  \sup_{x\in\R^d}
  \int_{K_{\ell}^c}\hspace{-0.2cm}
   R_1^{-V}(x,y)(R_1^{-V}(x,y)\land n)^{p-1}
  \m(\d y)
 \right|
\leq
 \sup_{x\in\R^d}
 \int_{|x-y|<\alpha_n}\hspace{-0.5cm}
  R_1(x,y)^{\,p}
 \m(\d y).
\label{eq:GreenEstimateSchroedinger}
\end{align}
The right-hand side of \eqref{eq:GreenEstimateSchroedinger}
converges to $0$ as $n\to\infty$ as shown above.
So we can obtain
$
 \m
\in
 S_{K_{\infty}}^{\,p}({\bf X}^{-V-1})
=
 S_{C\!K_{\infty}}^{\,p}({\bf X}^{-V-1})
$
from
$
 \m
\in
 S_{K_{\infty}}^1({\bf X}^{-V-1})
=
 S_{C\!K_{\infty}}^1({\bf X}^{-V-1})
$.
Therefore,
the assertion holds from Theorem~\ref{thm:compactEmbedding}.
\end{apf}
}
\end{ex}

\begin{ex}[Symmetric Relativistic \boldmath$\alpha$-stable Process]\label{ex:relastable}
{\rm
Take $0 <\alpha <2$ and $m\geq0$.
Let ${\bf X}=(\Omega,X_t,{\PP}_x)$
be a L\'evy process on $\R^d$ with
\begin{equation*}
 {\EE}_{0}
 \left[
  e^{\sqrt{-1}\<  \xi, X_{t}  \>}
 \right]
=
 \exp
 \left(
  -t\left\{(|\xi|^2+m^{2/\alpha})^{\alpha/2}-m\right\}
 \right).
\end{equation*}
If $m>0$, it is called the
\emph{relativistic $\alpha$-stable process with mass $m$}
(see \cite{CKS:POTA2012}).
In particular, if $\alpha=1$ and $m>0$, it is called the
\emph{free relativistic Hamiltonian process}
(see \cite{CMS:1990, CKS:AOP2012, HS:PP}).
When $m=0$, ${\bf X}$ is nothing but the usual
\emph{(rotationally) symmetric $\alpha$-stable process}.
It is known that ${\bf X}$ is transient if and only if
$d>2$ under $m>0$ or $d>\alpha$ under $m=0$, and
${\bf X}$ is a doubly Feller conservative process.

Let $(\E,\F)$ be the Dirichlet form on $L^2(\R^d)$
associated with ${\bf X}$.
Using Fourier transform
$
 \hat{f}(x)
:=
 \frac{1}{(2\pi)^{d/2}}\int_{\R^d}e^{i\<x,y\>}f(y){\rm d}y
$,
it follows from \cite[Example~1.4.1]{FOT} that
\begin{align*}
 \left\{
 \begin{array}{rl}
  \F
 \hspace{-3mm}&
 =
  \displaystyle{
  \left\{
   f\in L^2(\R^d)
  \;\left|\;
   \int_{\R^d}
    |\hat{f}(\xi)|^2
    \left((|\xi|^2+m^{2/\alpha})^{\alpha/2} -m \right)
   \d \xi
  <
   +\infty
  \right.
  \right\}
  }, \\
  \E(f,g)
 \hspace{-3mm}&
 =
  \displaystyle{
  \int_{\R^d}
   \hat{f}(\xi)\overline{\hat{g}(\xi)}
   \left((|\xi|^2+m^{2/\alpha})^{\alpha/2} -m \right)
  \d \xi
  \quad\text{ for }f,g\in \F.
  }
 \end{array}
 \right.
\end{align*}
Since ${\bf X}$ is a L\'evy process,
in view of \cite[Corollary~7.16]{BF:Potential},
there exists $C=C(m,\alpha)>0$ such that
\begin{align}
 \left\{
 \begin{array}{cl}
 \!\!\!\!\!\! &H^1(\R^d)\subset \F,   \\
 \!\!\!\!\!\! &\E_1(f,f)\leq C(\|\nabla f\|_2^2+\|f\|_2^2)\quad\text{ for }\quad f\in H^1(\R^d).
 \end{array}
 \right.
\label{eq:RelSobolevDir}
\end{align}
It is shown in \cite{CS3} that
the corresponding jumping measure $J$ of $(\E,\F)$  satisfies
\begin{align*}
 J(\d x \d y)
=
 J_m(x,y)\d x \d y
\quad \text{ with }~~J_m(x,y)
=
 A(d,-\alpha)
 \frac{\Psi(m^{1/\alpha}|x-y|)}{|x-y|^{d+\alpha}},
\end{align*}
where
$
 A(d,-\alpha)
=
 \frac{\alpha 2^{d+\alpha}\Gamma(\frac{d+\alpha}{2})}{2^{d+1}\pi^{d/2}\Gamma(1-\frac{\alpha}{2})}
$,
and $\Psi(r) := I(r)/I(0)$ with
\begin{equation*}
 I(r)
:=
 \int_0^{\infty}
  s^{\frac{d+\alpha}{2}-1}
  e^{-\frac{s}{4}-\frac{r^2}{s}}
 \d s
\end{equation*}
is a decreasing function satisfying
$\Psi(r)\asymp e^{-r}(1+r^{(d+\alpha-1)/2})$ near $r=+\infty$,
and $\Psi(r)=1+\Psi''(0)r^2/2+o(r^4)$ near $r=0$.
In particular,
\begin{align*}
 \left\{
 \begin{array}{rl}
  \F
 \hspace{-3mm}&
 =
  \displaystyle{
  \left\{
   f\in L^2(\R^d)
  \;\left|\;
   \int_{\R^d\times\R^d}|f(x)-f(y)|^2J_m(x,y) \d x \d y
  <
   +\infty
  \right.
  \right\}
  }, \\
  \E(f,g)
 \hspace{-3mm}&
 =
  \displaystyle{
  \frac{1}{2}
  \int_{\R^d\times\R^d}
   (f(x)-f(y))(g(x)-g(y))J_m(x,y)
  \d x \d y
  \quad\text{ for }\quad f,g\in \F
  }.
 \end{array}
 \right.
\end{align*}
Let $p_t(x,y)$ be the heat kernel of ${\bf X}$.
The following global heat kernel estimates are proved in
\cite[Theorem~2.1]{CKS:JFA2012}:
There exist $C_1,C_2>0$ such that
\begin{align}
 C_2^{-1}
 \Phi^m_{1/C_1}(t,x,y)
\leq
 p_t(x,y)
\leq
 C_2
 \Phi^m_{C_1}(t,x,y)
\quad\text{ for \ \  all }\quad (t,x,y)\in ]0,+\infty[\times\R^d\times\R^d,
\label{eq:GHKERS}
\end{align}
where
\begin{align*}
 \Phi^m_C(t,x,y)
:=
 \left\{
 \begin{array}{ll}
  t^{-d/\alpha} \wedge  tJ_m(x, y), & t\in]0,1/m], \\
  m^{d/\alpha-d/2} t^{-d/2}
  \exp
  \left(
   -C^{-1}(m^{1/\alpha}|x - y|
  \wedge
   m^{2/\alpha-1}\frac{|x-y|^2}{t})
  \right),
  & t\in]1/m,+\infty[.
 \end{array}
 \right.
\end{align*}
In particular, we have for $m=0$
\begin{align}
 p_t(x,y)
&\leq
 C_2
 t^{-d/\alpha}
\quad\text{ for \ \  all }\quad (t,x,y)\in]0,+\infty[\times\R^d\times\R^d,
\label{eq:UCHK}
\end{align}
and for $m>0$
\begin{align}
 p_t(x,y)
&\leq
 C_2
 m^{\frac{d}{\alpha}-\frac{d}{2}}
 (t^{-\frac{d}{\alpha}}+ t^{-\frac{d}{2}})
\quad\text{ for \ \  all }\quad
 (t,x,y)\in]0,+\infty[\times\R^d\times\R^d.
\label{eq:UCHKm}
\end{align}
It is shown in
\cite[Theorem~1.2 and Example~2.4]{CK:MixedHeat}
or
\cite[Theorem~1.2]{CKK:TranAMS2011,CKK:TranAMS2015} that
$p_t(x,y)$ is jointly continuous in
$(t,x,y)\in ]0,+\infty[\times \R^d\times \R^d$.
The $\beta$-order resolvent kernel
$
 R_{\beta}(x,y)
:=
 \int_0^{\infty}e^{-\beta t}p_t(x,y) \d t
\in
 [0,+\infty]
$
is also continuous in $(x,y)\in \R^d\times\R^d$.

We say that $\mu\in K_{d,\alpha}^{\,p}$ if and only if
\begin{align*}
 \lim_{r\to0}
 \sup_{x\in\R^d}
 \int_{|x-y|<r}\frac{\mu(\d y)}{|x-y|^{(d-\alpha)p}}=0
\quad &\text{ for }\quad d>\alpha,\\
 \lim_{r\to0}
 \sup_{x\in\R^d}
 \int_{|x-y|<r}(\log|x-y|^{-1})^p\mu(\d y)=0
\quad &\text{ for }\quad d=\alpha=1, \\
 \sup_{x\in\R^d}
 \int_{|x-y|\leq1}\mu(\d y)<+\infty
\quad &\text{ for }\quad d=1<\alpha.
\end{align*}
We write $K_{d,\alpha}$ instead of $K_{d,\alpha}^1$ for $p=1$.
Then we have $K_{d,\alpha}^{\,p}=S_K^{\,p}({\bf X})$ by
\cite[Theorem~3.1]{KwM:pKato}.
The $d$-dimensional Lebesgue measure $\m$ belongs to
$K_{d,\alpha}^{\,p}=S_K^{\,p}({\bf X})$
if and only if $p\in [1, d/(d-\alpha)_+[$
by \cite[Theorem~3.2 or Corollary~4.4]{KwM:pKato},
and for any non-negative bounded $g\in L^1(\R^d)$
the finite measure $g\m$ also belongs to
$S_{C\!K_{\infty}}^{\,p}({\bf X}^{(1)})$
(to $S_{C\!K_{\infty}}^{\,p}({\bf X})$
if ${\bf X}$ is transient) by Theorem~\ref{thm:Equals}.
Here
$d/(d-\alpha)_+ := d/(d-\alpha)$ if $d>\alpha$ and
$d/(d-\alpha)_+ := +\infty$ if $d\leq\alpha$.
The surface measure $\sigma_R$
on the $R$-sphere $\partial B_R(0)$ satisfies that
$\sigma_R(B_r(x))\leq C_2r^{d-1}$
for any $x\in\R^d$ and $r>0$ with some $C_2>0$, and
$\sigma_R(B_r(x))\geq C_1r^{d-1}$
for any $x\in\partial B_R(0)$ and $r\in]0,r_0[$
with some $C_1,r_0>0$.
By Theorem~\ref{thm:Equals} and
\cite[Theorems~3.1 and 3.2]{KwM:pKato},
we can conclude that
$\sigma_R\in S_{C\!K_{\infty}}^{\,p}({\bf X}^{(1)})$ holds
if and only if $p\in [1, (d-1)/(d-\alpha)_+[$,
moreover,
$\sigma_R\in S_{C\!K_{\infty}}^{\,p}({\bf X})$ holds
if and only if $p\in [1, (d-1)/(d-\alpha)_+[$
provided ${\bf X}$ is transient.
Here
$(d-1)/(d-\alpha)_+ := (d-1)/(d-\alpha)$ if $d>\alpha$ and
$(d-1)/(d-\alpha)_+ := +\infty$ if $d\leq\alpha$.

We consider a connected non-empty open set $D$ of $\R^d$.
The notion of regular point in $\partial D$
is similarly defined as in Example~\ref{ex:Brownian}.
Denote by $(\partial D)_r$
the set of regular points in boundary.
$D$ is said to be \emph{regular} if
$(\partial D)_r=\partial D$.
The part process  ${\bf X}_D=(\Omega, X_t^D,{\PP}_x)$
of ${\bf X}$ on $D$ is defined as the process
killed upon leaving $D$.
Let $R^D(x,y)$ be the Green function
with respect to ${\bf X}_D$.
$D$ is said to be \emph{Green-bounded} if
$
 \sup_{x\in D}
 \int_D R^D(x,y) \d y
=
 \sup_{x\in D}
 {\EE}_x[\tau_D]
<
 +\infty
$,
equivalently $\m\in S_{D_0}^1({\bf X}_D)$,
where $\m$ is the $d$-dimensional Lebesgue measure on $D$.
By \cite[Lemma~4.1]{KimKuwae:GeneralAnal},
$\m(D)<+\infty$ implies the Green-boundedness of $D$.
For a (positive) Radon measure $\nu$ on $\R^d$,
we consider the following conditions:
\begin{align}
 \text{ for }m=0,\quad
 \left\{
 \begin{array}{ll}
  d=1<\alpha, & D\text{ is bounded and }\nu(D)<+\infty, \\
  d=\alpha=1, & D\text{ is Green-bounded and } \nu\in S_K^{\,p}({\bf X})\text{ with }\nu(D)<+\infty, \\
  d>\alpha,   & \nu\in S_K^{\,p}({\bf X})\text{ with }\nu(D)<+\infty,
 \end{array}
 \right.
\label{eq:KatoGreenStable}
\end{align}
and
\begin{align}
 \text{ for }m>0,\quad
 \left\{
 \begin{array}{ll}
  d=1<\alpha, & D\text{ is bounded and }\nu(D)<+\infty, \\
  d=1\geq\alpha, \text{ or }d=2, & D\text{ is Green-bounded and } \nu\in S_K^{\,p}({\bf X})\text{ with }\nu(D)<+\infty, \\
  d\geq3, & \nu\in S_K^{\,p}({\bf X})\text{ with }\nu(D)<+\infty.
 \end{array}
 \right.
\label{eq:KatoGreenStable*}
\end{align}
As we see the following,
the difference of
\eqref{eq:KatoGreenStable} and
\eqref{eq:KatoGreenStable*} comes from
the order of $t$ in the upper heat kernels
\eqref{eq:UCHK} and \eqref{eq:UCHKm}.
It is unclear if $(x,y)\mapsto R^D(x,y)$
is extended continuous as in \cite[Theorem~2.10(iii)]{CZ:BS}.
So we do not know if
$(R^D)^p\nu\in C_b(D)$
($\in C_{\infty}(D)$ under the regularity of $D$)
under \eqref{eq:KatoGreenStable} or
\eqref{eq:KatoGreenStable*}.
However,
we can deduce the following:
\begin{prop}\label{prop:Coincidence}
Suppose that $\nu$ satisfies
\eqref{eq:KatoGreenStable} and \eqref{eq:KatoGreenStable*}.
Then
$
 \1_D\nu
\in
 S_{K_{\infty}}^{\,p}({\bf X}_D)
=
 S_{C\!K_{\infty}}^{\,p}({\bf X}_D)
$.
\end{prop}
\begin{pf}
Since
${\bf X}_D$ possesses {\bf (RSF)}, it suffices to show
$\1_D\nu\in S_{K_{\infty}}^{\,p}({\bf X}_D)$
by Theorem~\ref{thm:Coincidence}.
First
we prove $\1_D\nu\in S_K^{\,p}({\bf X}_D)$ under
\eqref{eq:KatoGreenStable} and \eqref{eq:KatoGreenStable*}.
Consider the case $d=1<\alpha$ in both cases.
Since $\nu(D)<+\infty$, we see
$\1_D\nu\in K_{1,\alpha}^{\,p}=S_K^{\,p}({\bf X})$
so that $\1_D\nu\in S_K^{\,p}({\bf X}_D)$.
In other cases,
$\1_D\nu\in S_K^{\,p}({\bf X}_D)$ follows from
$\nu\in S_K^{\,p}({\bf X})$.
So
it suffices to show the $L^p$-Green-tightness of $\1_D\nu$
under ${\bf X}_D$ in the sense of Zhao.
To do this, we prove the $L^p$-Green-tightness of $\1_D\nu$
under ${\bf X}_D^{(1)}$ in the sense of Zhao.
By use of the claims (C1)--(C4) in
\cite[the proof of Theorem~4.1]{KimKuwae:GeneralAnal},
we have that
\begin{align}
 R^D_{\beta}(x,y)
&=
 R^D_{\beta}(x,y)\land n
\quad\text{ for }\quad |x-y|\geq \alpha_n,
\label{eq:CutOff}\\
 R^D(x,y)
&=
 R^D(x,y)\land n
\quad\text{ for }\quad |x-y|\geq \beta_n,
\label{eq:CutOff0}
\end{align}
where
$
 \alpha_n
:=
 \left(\frac{C_1}{n-C_2}\right)^{\frac{1}{d-\alpha}}
$
for $d>\alpha$, $n>C_2$,
$
 \alpha_n
:=
 \exp\left(-\frac{n-C_2}{C_1}\right)
$
for $d=\alpha=1$, $n>C_2$,
and $\alpha_n>0$ is arbitrary for $d=1<\alpha$, and
$\beta_n := \left(\frac{C_3}{n}\right)^{\frac{1}{d-\alpha}}$
for $d>\alpha$ with $m=0$
and
$
 \beta_n
:=
 \left(
  \frac{2C_3}{n}
 \right)^{\frac{1}{d-\alpha}}
\lor
 \left(
  \frac{2C_3}{n}
  m^{\frac{2-\alpha}{\alpha}}
 \right)^{\frac{1}{d-2}}
$
for $d\geq3$ with $m>0$.
Here $C_1, C_2$ and $C_3$ are the positive constants
appeared in the claims (C1)--(C4) in
\cite[the proof of Theorem~4.1]{KimKuwae:GeneralAnal},
$C_2$ depends on $\beta$ and $C_3$ depends on $m$.
(Note that there is a typo in (C4);
the upper bound of $R(x, y)$ for $m>0$ and $d\geq 3$
should be
$
 C_3
 (1+m^{\frac{2-\alpha}{\alpha}} |x-y|^{2-\alpha})
 /|x-y|^{d-\alpha}
$.)
From \eqref{eq:CutOff}, for any $A\in\mathscr{B}(\R^d)$,
we have
\begin{align}
 \left|
  \sup_{x\in D}
  \int_A R_{\beta}^D(x,y)^{\,p}\nu(\d y)
 -
  \sup_{x\in D}
  \int_A(R_{\beta}^D(x,y)\land n)^{\,p}\nu(\d y)
 \right|
\leq
 \sup_{x\in D}
 \int_{\{|x-y|<\alpha_n\}}
  R_{\beta}(x,y)^{\,p}
 \nu(\d y).
\label{eq:GreenEst2}
\end{align}
Since $\nu\in S_K^{\,p}({\bf X})=K_{d,\alpha}^{\,p}$,
the right-hand side of \eqref{eq:GreenEst2}
converges to $0$ by \cite[Theorem~3.1]{KwM:pKato}.
Taking an increasing sequence $\{K_{\ell}\}$ of
compact subsets of $D$, we have
\begin{align*}
 \lim_{\ell\to\infty}
 \sup_{x\in D}&
 \int_{D\setminus K_{\ell}}
  R^D_{\beta}(x,y)^{\,p}
 \nu(\d y)\\
&\leq
 \lim_{\ell\to\infty}
 \sup_{x\in D}
 \int_{D\setminus K_{\ell}}
  (R^D_{\beta}(x,y)\land n)^{\,p}
 \nu(\d y)
+
 \sup_{x\in D}
 \int_{\{|x-y|<\alpha_n\}}
  R_{\beta}(x,y)^{\,p}
 \nu(\d y)
\nonumber\\
&\leq
 \sup_{x\in \R^d}
 \int_{\{|x-y|<\alpha_n\}}
  R_{\beta}(x,y)^{\,p}
 \nu(\d y)
\longrightarrow
 0
\quad\text{ as }\quad n\to\infty.
\end{align*}
Thus $\1_D\nu\in S_{K_{\infty}}^{\,p}({\bf X}_D^{(1)})$.
If $D$ is Green-bounded, i.e. $\1_D\m\in S_{D_0}({\bf X}_D)$,
then we obtain $\1_D\nu\in S_{K_{\infty}}^{\,p}({\bf X}_D)$
from Proposition~\ref{prop:GreenBddEqui}.
From now on,
we consider the transient case, i.e.,
$d>\alpha$ with $m=0$ or $d\geq3$ with $m>0$.
In this case,
by replacing \eqref{eq:CutOff} with \eqref{eq:CutOff0},
a similar estimate with \eqref{eq:GreenEst2} holds
in the following manner: for any $A\in\mathscr{B}(\R^d)$
\begin{align}
 \left|
  \sup_{x\in D}
  \int_A R^D(x,y)^{\,p}\nu(\d y)
 -
  \sup_{x\in D}
  \int_A(R^D(x,y)\land n)^{\,p}\nu(\d y)
 \right|
\leq
 \sup_{x\in D}
 \int_{\{|x-y|<\beta_n\}}
  R(x,y)^{\,p}
 \nu(\d y).
\label{eq:GreenEst0}
\end{align}
The right-hand side of \eqref{eq:GreenEst0} converges to $0$,
because $\nu\in S_{K}^{\,p}({\bf X})=K_{d,\alpha}^{\,p}$
and the fourth claim (C4) yields
$R(x,y)\leq C/|x-y|^{d-\alpha}=C G(x, y)$ for $|x-y|<\beta_n$
with small $\beta_n$.
Taking an increasing sequence of compact sets as above
\begin{align*}
 \lim_{\ell\to\infty}
 \sup_{x\in D}&
 \int_{D\setminus K_{\ell}}
  R^D(x,y)^{\,p}
 \nu(\d y)\\
&\leq
 \lim_{\ell\to\infty}
 \sup_{x\in D}
 \int_{D\setminus K_{\ell}}
  (R^D(x,y)\land n)^{\,p}
 \nu(\d y)
+
 \sup_{x\in D}
 \int_{\{|x-y|<\alpha_n\}}
  R(x,y)^{\,p}
 \nu(\d y)
\nonumber\\
&\leq
 \sup_{x\in \R^d}
 \int_{\{|x-y|<\alpha_n\}}
  R(x,y)^{\,p}
 \nu(\d y)
\longrightarrow
 0
\quad\text{ as }\quad n\to\infty.
\end{align*}
Therefore
$
 \1_D\nu
\in
 S_{K_{\infty}}^{\,p}({\bf X}_D)
=
 S_{C\!K_{\infty}}^{\,p}({\bf X}_D)
$.
\end{pf}
It is proved in \cite[Corollary~3.3]{KwM:pKato} that
$p(d-\alpha)<d$ is equivalent to $\m\in S_K^{\,p}({\bf X})$.
From this, if
\begin{align}
\text{for } m=0,\quad
 \left\{
 \begin{array}{ll}
  d=1\leq\alpha, & D\text{ is bounded}, \\
  d>\alpha, & p(d-\alpha)<d \text{ with }\m(D)<+\infty,
 \end{array}
 \right.
\label{eq:KatoGreenUnderlyingStable}
\end{align}
and
\begin{align}
 \text{for }m>0,\quad
 \left\{
 \begin{array}{ll}
  d=1,    & p(d-\alpha)<d \text{ and }D\text{ is bounded}, \\
  d\geq2, & p(d-\alpha)<d \text{ with }\m(D)<+\infty,
 \end{array}
 \right.
\label{eq:KatoGreenUnderlyingStable*}
\end{align}
then
$
 \m
\in
 S_{K_{\infty}}^{\,p}({\bf X}_D)
=
 S_{C\!K_{\infty}}^{\,p}({\bf X}_D)
$
provided $D$ is a regular domain.
We relax \eqref{eq:KatoGreenUnderlyingStable} and
\eqref{eq:KatoGreenUnderlyingStable*} in the following:
\begin{align}
 \text{for }m=0,\quad
 \left\{
 \begin{array}{ll}
  d=1\leq\alpha, & D\text{ is bounded}, \\
  d>\alpha,      & p(d-\alpha)<d \text{ with }\lim_{x\in D,|x|\to\infty}\m(D\cap B_1(x))=0,
 \end{array}
 \right.
\label{eq:KatoGreenUnderlyingStableRrelaxed}
\end{align}
and
\begin{align}
 \text{for }m>0,\quad
 \left\{
 \begin{array}{ll}
  d=1,    & p(d-\alpha)<d\text{ and }D\text{ is bounded}, \\
  d\geq2, & p(d-\alpha)<d \text{ with }\lim_{x\in D,|x|\to\infty}\m(D\cap B_1(x))=0.
 \end{array}
 \right.
\label{eq:KatoGreenUnderlyingStableRelaxed*}
\end{align}
\begin{lem}\label{lem:ClassT}
Suppose that $\lim_{x\in D,|x|\to\infty}\m(D\cap B_1(x))=0$.
Then for the symmetric relativistic $\alpha$-stable process ${\bf X}$
the absorbing process ${\bf X}_D$ of ${\bf X}$
killed upon leaving $D$ is in Class {\bf (T)}
defined in Section \ref{sec:introduction}.
In particular,
$
 \1_D\m
\in
 S_{K_{\infty}}^1({\bf X}_D^{(1)})
=
 S_{C\!K_{\infty}}^1({\bf X}_D^{(1)})
$.
\end{lem}
\begin{pf}
The proof is similar with the proofs of
\cite[Lemmas~3.2 and 3.3]{TakTawTsuchi:Compact}
done for Brownian motion.
Since ${\bf X}$ is  ultra-contractive for small time,
i.e. $p_t(x,y)\leq C/t^{d/\alpha}$ holds
for all $x,y\in \R^d$ and $t\in]0,1/m[$,
the proofs of
\cite[Lemmas~3.2 and 3.3]{TakTawTsuchi:Compact}
remain valid
for symmetric relativistic $\alpha$-stable processes,
i.e., we obtain $\lim_{x\in D,|x|\to\infty}P_t^D1(x)=0$
for $t\in]0,1/m[$, and this also holds
for general $t\in]0,+\infty[$.
Note that
the proof of \cite[Lemma~3.2]{TakTawTsuchi:Compact}
relies on the translation invariance of ${\bf X}$.
\end{pf}
Therefore,
by Theorem~\ref{thm:compactEmbedding} we have the following:
\begin{thm}\label{thm:CptEmbedStableEuclidean}
Suppose that
\eqref{eq:KatoGreenUnderlyingStableRrelaxed} and
\eqref{eq:KatoGreenUnderlyingStableRelaxed*} are satisfied.
Then
$
 \1_D\m
\in
 S_{K_{\infty}}^{\,p}({\bf X}_D^{(1)})
=
 S_{C\!K_{\infty}}^{\,p}({\bf X}_D^{(1)})
$,
in particular,  the embedding
\begin{align*}
 \F_D\hookrightarrow L^{2p}(D;\m)
\end{align*}
is compact.
Moreover, if $D$ is Green-bounded, then
$
 \1_D\m
\in
 S_{K_{\infty}}^{\,p}({\bf X}_D)
=
 S_{C\!K_{\infty}}^{\,p}({\bf X}_D)
$,
in particular, the embedding
\begin{align*}
 (\F_D)_e\hookrightarrow L^{2p}(D;\m)
\end{align*}
is compact.
Here $(\F_D)_e$ is the extended Dirichlet space of
$(\E,\F_D)$ on $L^2(D)$.
\end{thm}

\begin{pf}
The latter assertion follows from
Proposition~\ref{prop:GreenBddEqui}.
So it suffices to prove the former assertion.
The condition $p(d-\alpha)<d$ in
\eqref{eq:KatoGreenUnderlyingStableRrelaxed} and
\eqref{eq:KatoGreenUnderlyingStableRelaxed*} is
equivalent to $\m\in S_K^{\,p}({\bf X})$,
hence it implies $\1_D\m\in S_K^{\,p}({\bf X}_D)$.
By using Lemma~\ref{lem:ClassT},
the proof is similar to the proof of
Theorem~\ref{thm:CptEmbedEuclideanRelaxed}.
To do this, from claims (C1)--(C3) in
\cite[proof of Theorem~4.1]{KimKuwae:GeneralAnal} and
\eqref{eq:GHKERS},
we only note the following estimate;
there exists $C>0$ which depends on $d$, $\alpha$ and $m$
such that
\begin{align*}
 R_1^D(x,y)
\leq
 R_1(x,y)
\leq
 \left\{
 \begin{array}{lr}
  C,
  & d=1<\alpha, \\
  \displaystyle
  C\left(\frac{1}{|x-y|^2}+1\right),
  & d=\alpha=1, \\
  \displaystyle
  C\left(\frac{1}{|x-y|^{d-\alpha}}+1\right),
  & d>\alpha
 \end{array}
 \right.
\end{align*}
for all $x, y\in D$.
\end{pf}

Recall $\mathscr{B}_0$ defined in \eqref{eq:B0}.
As proved in Lemma~\ref{lem:GrennTightBM}, we can prove that
$B\in\mathscr{B}_0$ implies
$
 \m^B
:=
 \1_B\m
\in
 S_{K_{\infty}}^1({\bf X}^{(1)})
=
 S_{C\!K_{\infty}}^1({\bf X}^{(1)})
$
by using \eqref{eq:UCHK} and \eqref{eq:UCHKm}.
We now claim the following:
\begin{prop}\label{prop:compactEmbedEucStable}
Suppose $p(d-\alpha)<d$ and $B\in\mathscr{B}_0$.
Then
$
 \m^B
:=
 \1_B\m
\in
 S_{K_{\infty}}^{\,p}({\bf X}^{(1)})
=
 S_{C\!K_{\infty}}^{\,p}({\bf X}^{(1)})
$.
In particular,
$\F$ is compactly embedded into $L^{2p}(B)$.
\end{prop}
\begin{pf}
It suffices to prove
$
 \m^B
:=
 \1_B\m
\in
 S_{K_{\infty}}^{\,p}({\bf X}^{(1)})
=
 S_{C\!K_{\infty}}^{\,p}({\bf X}^{(1)})
$
by
Theorems~\ref{thm:compactEmbedding} and \ref{thm:Coincidence}.
By \eqref{eq:GreenEstimates}, we have
$R_1(x,y)=R_1(x,y)\land n$ for $|x-y|\geq\alpha_n$,
where $\alpha_n$ is a sequence converging to $0$,
which can be constructed based on the estimates
in the proof of Theorem~\ref{thm:CptEmbedStableEuclidean}.
Then one can deduce the following estimate:
\begin{align}
 \left|
  \sup_{x\in\R^d}
  \int_{K_{\ell}^c}\hspace{-0.2cm}
   R_1(x,y)^{\,p}
  \m^B(\d y)
 -
  \sup_{x\in\R^d}
  \int_{K_{\ell}^c}\hspace{-0.2cm}
   R_1(x,y)(R_1(x,y)\land n)^{p-1}
  \m^B(\d y)
 \right|
\leq
 \sup_{x\in\R^d}
 \int_{|x-y|<\alpha_n}\hspace{-0.5cm}
  R_1(x,y)^{\,p}
 \m(\d y).
\label{eq:GreenEstimateStable}
\end{align}
The right-hand side of \eqref{eq:GreenEstimateBM}
converges to $0$ as $n\to\infty$
from $\m\in S_K^{\,p}({\bf X})$ under $p(d-\alpha)<d$
by applying \cite[Theorem~3.1]{KwM:pKato} except $d=1<\alpha$.
When $d=1<\alpha$,
the right-hand side of \eqref{eq:GreenEstimateStable}
is estimated above by
$C(m)^p\m(B_{\alpha_n}(0))$, which goes to $0$ as $n\to\infty$.
Here $C(m) := \sup_{x,y\in\R^d}R_1(x,y)>0$ is the constant
for the case $d=1<\alpha$.
Thus we can obtain the assertion.
\end{pf}
\begin{cor}\label{cor:EdmundEvanceStable}
Let $D$ be a domain of $\R^d$.
Suppose $p(d-\alpha)<d$.
Then the following statements are equivalent:
\begin{enumerate}
\item
$D\in\mathscr{B}_0$,\quad
\item
$
 \1_D\m
\in
 S_{K_{\infty}}^{\,p}({\bf X}^{(1)})
=
 S_{C\!K_{\infty}}^{\,p}({\bf X}^{(1)})
$,\quad
\item
$\F$ is compactly embedded into $L^{2p}(D)$.
\end{enumerate}
\end{cor}
\begin{pf}
The condition $p(d-\alpha)<d$ is only used
to establish the continuity of the embedding
$\F\hookrightarrow L^{2p}(D)$.
We have already proved
(1)$\Longrightarrow$(2)$\Longrightarrow$(3)
in Proposition~\ref{prop:compactEmbedEucStable}.
The proof of the implication
(3)$\Longrightarrow$(1) is similar to
\cite[Chapter X, Lemma~6.11]{EdEv:Spectral}
by using \eqref{eq:RelSobolevDir}.
\end{pf}

In the end of this example,
we give the compactness of the Schr\"odinger semigroups,
which is a $p$-version of
\cite[Theorems~5.2]{TakTawTsuchi:Compact}
to symmetric relativilistic $\alpha$-stable processes.
\begin{thm}\label{thm:SchoedingerCopactStable}
Let $V$ be a positive Borel function on $\R^d$
satisfying $V\m\in S_{L\!K}^1({\bf X})$.
Suppose that $\{V\leq M\}\in\mathscr{B}_0$
for any $M>0$ and $p(d-\alpha)<d$.
Then the Schr\"odinger semigroup $P_t^{-V}$
defined in \eqref{eq:SchroedingerSemigroup}
forms a compact operator from $L^2(\R^d)$ to $L^{2p}(\R^d)$.
\end{thm}

\begin{cor}\label{cor:SchoedingerCopactStable}
Let $V$ be a positive continuous function on $\R^d$
satisfying $V\m\in S_{L\!K}^1({\bf X})$.
Suppose $p(d-\alpha)<d$ and
\begin{align*}
 \lim_{|x|\to\infty}V(x)=+\infty.
\end{align*}
Then the Schr\"odinger semigroup $P_t^{-V}$
defined in \eqref{eq:SchroedingerSemigroup}
forms a compact operator from $L^2(\R^d)$ to $L^{2p}(\R^d)$.
\end{cor}

\begin{rem}\label{rem:Matsuura2}
{\rm
Theorem~\ref{thm:SchoedingerCopact} and
Corollary~\ref{cor:SchoedingerCopactStable}
are not included in \cite[Example~4.4]{MatsuKohei:compact},
because we consider relativistic $\alpha$-stable process
${\bf X}$,
and the compactness of $P_t^{-V}$ is a different type from in
\cite[Theorems~2.2 and 2.4]{MatsuKohei:compact}.
}
\end{rem}

\begin{apf}{Theorem~\ref{thm:SchoedingerCopactStable}}
Let ${\bf X}^{-V-1}$ be the subprocess of ${\bf X}$
by $\exp\left(-\int_0^tV(X_s)\d s-t \right)$,
which possesses {\bf (RSF)} because of
$V\in S_{L\!K}^1({\bf X})$.
Though the framework of
\cite[Theorems~5.1 and 5.2]{TakTawTsuchi:Compact}
treats only the case $m=0$,
the proof remains valid for $m\geq 0$
in view of \eqref{eq:UCHK} and \eqref{eq:UCHKm}.
Then
the assertion for $p=1$ is nothing but
\cite[Theorem~5.2]{TakTawTsuchi:Compact},
which was done by showing
$
 \m
\in
 S_{K_{\infty}}^1({\bf X}^{-V-1})
=
 S_{C\!K_{\infty}}^1({\bf X}^{-V-1})
$.
Suppose $p>1$ and
let $R_1^{-V}(x,y)$ be the Green kernel of ${\bf X}^{-V-1}$.
Then
\begin{align}
 \left|
  \sup_{x\in\R^d}
  \int_{K_{\ell}^c}\hspace{-0.2cm}
   R_1^{-V}(x,y)^{\,p}
  \m(\d y)
 -
  \sup_{x\in\R^d}
  \int_{K_{\ell}^c}\hspace{-0.2cm}
   R_1^{-V}(x,y)(R_1^{-V}(x,y)\land n)^{p-1}
  \m(\d y)
 \right|
\leq
 \sup_{x\in\R^d}
 \int_{|x-y|<\alpha_n}\hspace{-0.5cm}
  R_1(x,y)^{\,p}
 \m(\d y).
\label{eq:GreenEstimateSchroedingerStable}
\end{align}
The right-hand side of
\eqref{eq:GreenEstimateSchroedingerStable}
converges to $0$ as $n\to\infty$ as shown above.
So we can obtain
$
 \m
\in
 S_{K_{\infty}}^{\,p}({\bf X}^{-V-1})
=
 S_{C\!K_{\infty}}^{\,p}({\bf X}^{-V-1})
$
from
$
 \m
\in
 S_{K_{\infty}}^1({\bf X}^{-V-1})
=
 S_{C\!K_{\infty}}^1({\bf X}^{-V-1})
$.
Therefore,
the assertion holds from Theorem~\ref{thm:compactEmbedding}.
\end{apf}
}
\end{ex}

\providecommand{\bysame}{\leavevmode\hbox to3em{\hrulefill}\thinspace}
\providecommand{\MR}{\relax\ifhmode\unskip\space\fi MR }
\providecommand{\MRhref}[2]{%
  \href{http://www.ams.org/mathscinet-getitem?mr=#1}{#2}
}
\providecommand{\href}[2]{#2}


\begin{thebibliography}{99}




\bibitem{BF:Potential}
  {C.~Berg and G.~Forst}, 
  \emph{Potential theory on locally compact abelian groups},   
  Ergebnisse der Mathematik und ihrer Grenzgebiete, Band 87. Springer-Verlag, New York-Heidelberg, 1975. 




\bibitem{Chen:gaugeability2002}
  {Z.-Q.~Chen}, 
  \emph{Gaugeability and conditional gaugeability}, 
  Trans. Amer. Math. Soc. {\bf 354} (2002), no. 11, 4639--4679.

\bibitem{CMS:1990}
  {R.~Carmona, W.~C.~Masters and B.~Simon},  
  \emph{Relativistic Schr\"odinger operators: Asymptotic behavior of the eigenfunctions},
  J. Funct. Anal. {\bf 91} (1990), no.~1, 117--142.

\bibitem{CKK:TranAMS2011}
  {Z.-Q.~Chen, P.~Kim and T.~Kumagai}, 
  \emph{Global heat kernel estimates for symmetric jump processes}, 
  Trans. Amer. Math. Soc. {\bf 363} (2011), no.~9, 5021--5055. 

\bibitem{CKK:TranAMS2015}
  \bysame, 
  \emph{Corrigendum to \lq\lq Global heat kernel estimates for symmetric jump processes'' [MR2806700]},
  Trans. Amer. Math. Soc. {\bf 367} (2015), no. 10, 7515.

\bibitem{CKS:AOP2012}
  {Z.-Q. Chen, P.~Kim and R.~Song}, 
  \emph{Sharp heat kernel estimates for relativistic stable processes in open sets},
  Ann. Prob. {\bf 40} (2012), no.~1, 213--244.

\bibitem{CKS:POTA2012}
  {Z.-Q. Chen, P.~Kim and R.~Song}, 
  \emph{Global heat kernel estimates for relativistic stable processes in half-space-like open sets},
  Potential Anal. {\bf 36} (2012), no. 2, 235--261.

\bibitem{CKS:JFA2012}
  \bysame, 
  \emph{Global heat kernel estimates for relativistic stable processes in exterior open sets}. 
  J. Funct. Anal. {\bf 263} (2012), no. 2, 448--475.

\bibitem{CK:MixedHeat}
  {Z.-Q.~Chen and T.~Kumagai},     
  \emph{Heat kernel estimates for jump processes of mixed types on metric measure spaces}, 
  Probab. Theory Related Fields {\bf 140} (2008), no. 1-2, 277--317.

\bibitem{CS3}
  {Z.-Q. Chen and R. Song},
  \emph{Drift transforms and Green function estimates for discontinuous processes}, 
  {J. Funct. Anal.} {\bf 201} (2003), 262--281.

\bibitem{CSgauge2002}    
  \bysame, 
  \emph{General gauge and conditional gauge theorems},
  Ann. Prob. {\bf 30} (2002), no.~3, 1313--1339. 

\bibitem{CSgauge2003}    
  \bysame, 
  \emph{Conditional gauge theorem for non-local Feynman-Kac transforms},
  Probab. Theory Related Fields {\bf 125} (2003), no. 1, 45--72.  



\bibitem{Chungdoubl} 
  {K.~L. Chung},
  \emph{Doubly-Feller process with multiplicative functional,}
  Seminar on stochastic processes, 1985 (Gainesville, Fla., 1985), 63--78, Progr. Probab. Statist. {\bf 12}, Birkh\"auser Boston, Boston, MA, 1986.

\bibitem{CZ:BS}
  {K. L. Chung and Z. Zhao},
  \emph{From Brownian motion to Schr\"odinger's equation}. 
  Grundlehren der Mathematischen Wissenschaften [Fundamental Principles of Mathematical Sciences], 312. Springer-Verlag, Berlin, 1995.




\bibitem{EdEv:Spectral}
  {D.~E.~Edmunds and W.~D.~Evans}, 
  \emph{Spectral theory and differential operators}, 
  Oxford Mathematical Monographs. Oxford Science Publications. 
  The Clarendon Press, Oxford University Press, New York, 1987.




\bibitem{FOT}
  {M. Fukushima, Y. Oshima and M. Takeda},
  \emph{Dirichlet forms and symmetric Markov processes}. 
  de Gruyter Studies in Mathematics, {\bf 19} 
  Walter de Gruyter \& Co., Berlin, 1994.

   


  





\bibitem{HS:PP}
  {I. W. Herbst and A. D. Sloan},
  \emph{Perturbation of translation invariant positivity preserving semigroups on $L^2(\R^n)$},    
  Trans. Amer. Math. Soc. {\bf 236} (1978), 325--360.






\bibitem{KK:AnalChara}
  {D. Kim and K. Kuwae},  
  \emph{Analytic characterizations of gaugeability for generalized Feynman-Kac functionals},
  Trans. Amer. Math. Soc. {\bf 369} (2017), no.~7, 4545--4596. 

\bibitem{KimKuwae:GeneralAnal}
  \bysame,
  \emph{General analytic characterization of gaugeability for Feynman-Kac functionals},
  Math. Ann. {\bf 370} (2018), no.~1-2, 1--37.

\bibitem{KKK:Kyoto}
  M. Kurniawaty, K. Kuwae and K. Tsuchida,
  \emph{On the doubly Feller property of resolvent},
  Kyoto J. Math. {\bf 57} (2017), no. 3, 637--654.

\bibitem{KwM:pKato}
  {K.~Kuwae and T.~Mori},
  \emph{$L^p$-Kato class measures for symmetric Markov processes under heat kernel estimates},
  2020, preprint.
  {\tt https://arxiv.org/pdf/2008.10934.pdf}

\bibitem{KwT:Katofunc}
  {K. Kuwae and M. Takahashi}, 
  \emph{Kato class functions of Markov processes under ultracontractivity},
  Potential Theory in Matsue, 193--202, \emph{Adv. Stud. Pure Math}., {\bf 44}, Math. Soc. Japan, Tokyo, 2006. 

\bibitem{KwT:Katounderheat}
  \bysame, 
  \emph{Kato class measures of symmetric Markov processes under heat kernel estimate},
  \emph{J. Funct. Anal.} {\bf 250} (2007), no. 1, 86--113.

 
\bibitem{MatsuKohei:compact} 
  {K.~Matsuura}, 
  \emph{Compactness of semigroups of explosive symmetric Markov processes},
  preprint, to appear in Kyoto J. Math. 2020. 

\bibitem{TM:pKato}
  {T.~Mori}, 
  \emph{$L^p$-Kato class measures and their relations with Sobolev embedding theorems},
  2020, preprint. 
  {\tt https://arxiv.org/pdf/2005.13758v2.pdf} 
 






\bibitem{Schilling:BM}
  {R,~Schilling}, 
  \emph{Measures, integrals and martingales},
  Cambridge University Press, New York, 2005.

\bibitem{SV;potential}
  {P.~Stollmann and J.~Voigt},
  \emph{Perturbation of Dirichlet forms by measures},
  Potential Anal. {\bf 5} (1996), no.~2, 109--138.   
 


	 




\bibitem{Takeda:Compact}
  {M.~Takeda},
  \emph{Compactness of symmetric Markov semi-groups and boundedness of eigenfunctions}, 
  Trans. Amer. Math. Soc. {\bf 372} (2019), no.~6, 3905--3920.

\bibitem{TakTawTsuchi:Compact}
  {M.~Takeda, Y.~Tawara and K.~Tsuchida},
  \emph{Compactness of Markov and Schr\"odinger semi-groups: A probabilistic approach},
  Osaka J. Math. {\bf 54} (2017), no.~3, 517--532.
   

\bibitem{Terkelsen:1972}  
  {F.~Terkelsen}, 
  \emph{Some minimax theorems}, 
  Math. Scand. {\bf 31} (1972), 405--413.


\bibitem{Var85}
  {N.~T.~Varopoulos}, 
  \emph{Hardy-Littelwood theory for semigroups},
  J. Func. Anal. {\bf 63} (1985), no.~2, 240--260.

\bibitem{Zhao:1983}
  {Z.~Zhao}, 
  \emph{Conditional gauge with unbounded potential}, 
  Z. Wahrscheinlichkeitstheorie verw. Gebiete
  {\bf 65} (1983), 13--18.   
  
\bibitem{Zhao:1989}
  \bysame,
  \emph{Gaugeability for unbounded domains}. 
  Seminar on Stochastic Processes (1989), Birkh\"auser, Boston, 207--214.

\bibitem{Zhao:Subcri}
  \bysame,
  \emph{Subcriticality and gaugeability of the Schr\"odinger operator},
  Trans. Amer. Math. Soc. {\bf 334} (1992), no. 1, 75--96.

\bibitem{Yan:1988}
  {J.~A. Yan},
  \emph{A formula for densities of transition functions},
  In S\'{e}minaire de {P}robabilit\'{e}s, {XXII}, volume 1321 of
  {\em Lecture Notes in Math.}, pages 92--100. Springer, Berlin, 1988.
\end{thebibliography}
\end{document}